\documentclass[]{article}
\usepackage{latexsym,diagram}

\newcommand{\kitem}{\begin{itemize}\vspace{-2ex}}
\newcommand{\kenditem}{\vspace{-1ex}\end{itemize}}
\newcommand{\bean}{\[\begin{array}{rcl}}
\newcommand{\eean}{\end{array}\]}

\newcommand{\Z}{{Z\!\!\!Z}}
\newcommand{\Q}{{I\!\!\!\!Q}}
\newcommand{\R}{{I\!\!R}}
\newcommand{\C}{{\,I\!\!\!\!C}}
\newcommand{\PP}{I\!\!P}
\newcommand{\CA}{{\cal A}}
\newcommand{\CB}{{\cal B}}
\newcommand{\CC}{{\cal C}}

\newcommand{\CH}{{\cal H}}
\newcommand{\CM}{{\cal M}}
\newcommand{\CF}{{\cal F}}

\newcommand{\Spec}{\mbox{\rm Spec}}
\renewcommand{\span}{\mbox{\rm span}}

\newcommand{\innt}{\mbox{\rm int}}
\newcommand{\cone}{\mbox{\rm cone}}

\newcommand{\surj}{\rightarrow\hspace{-0.8em}\rightarrow}

\newcommand{\ko}{\overline}
\newcommand{\kov}[3]{
   \unitlength=1.0pt \begin{picture}(0,0)(0,0)
   \put(0,#2){\vector(1,0){#1}}
   \end{picture}#3}
\newcommand{\ku}{\underline}
\newcommand{\ks}{\scriptstyle}
\newcommand{\kss}{\scriptscriptstyle}

\newcommand{\kf}{\footnotesize}
\newcommand{\keps}{\varepsilon}
\newcommand{\veee}{{\scriptscriptstyle\vee}}
\newcommand{\kb}{{\kss \bullet}}

\newcommand{\bHom}{\ko{\mbox{\rm Hom}}\,}

\reversemarginpar
\newcommand{\kkk}[1]{}
\newcommand{\ktrash}[1]{{}}

\newcommand{\AB}{{\CA\!\CB}}
\newcommand{\kk}{{K}}

\newcommand{\cx}{{\Sigma}}
\newcommand{\cxV}{{V}}
\newcommand{\cxD}{{d}}
\newcommand{\cxa}{{\sigma}}
\newcommand{\cxb}{{\tau}}
\newcommand{\cxF}{{\CF}}
\newcommand{\cxM}[1]{{M^{#1}}}
\newcommand{\Mi}{i}
\newcommand{\Mj}{j}
\newcommand{\Mk}{k}
\newcommand{\Mn}{N}
\newcommand{\nM}{{\CM}}

\newcommand{\cf}{{\cx}}
\newcommand{\cfV}{{\cxV}}
\newcommand{\cfD}{{\cxD}}
\newcommand{\cfa}{{\cxa}}
\newcommand{\cfb}{{\cxb}}
\newcommand{\cfF}{{\cxF}}
\newcommand{\cfM}[1]{{\cxM{#1}}}

\newcommand{\cxCF}[1]{{C^{#1}(\cx,\cxF)}}
\newcommand{\cxHF}[1]{{H^{#1}(\cx,\cxF)}}

\newcommand{\kP}{{\Delta}}
\newcommand{\kDP}{{\Diamond}}
\newcommand{\kT}{{D}}
\newcommand{\kTsh}{{\cal \kT}}
\newcommand{\kl}{{\ell}}
\newcommand{\kPM}[1]{{\cfM{#1}}}

\newcommand{\kn}{{k}}
\newcommand{\ki}{{p}}
\newcommand{\kj}{{q}}

\newcommand{\kPV}[1]{{V^{#1}}}
\newcommand{\Vi}{{k}}

\newcommand{\Vn}{{M}}
\newcommand{\kPP}[2]{{M^{#1}_{#2}}}
\newcommand{\kPPb}[2]{{\bar{M}^{#1}_{#2}}}
\newcommand{\ai}{{\nu}}
\newcommand{\aj}{{\mu}}

\newcommand{\an}{{m}}
\newcommand{\kPF}{{F}}
\newcommand{\kPFb}{{\bar{\kPF}}}

\newcommand{\kPX}{{X}}
\newcommand{\kPXb}{{\bar{\kPX}}}
\newcommand{\kPMb}[1]{{\bar{\cfM{}}^{#1}}}
\newcommand{\kts}[3]{{s(#1,#3,#2)}}

\newcommand{\fcP}{{\ku{\cone}(\kP)}}
\newcommand{\fcPd}{{\ku{\cone}(\kP^\veee)}}
\newcommand{\fsigma}{{\ku{\sigma}}}

\newcommand{\fn}{{\cf}}
\newcommand{\fnP}{{\fn(\kP)}}

\newcounter{Abschnitt}[section]
\newcommand{\neu}[1]{\protect\refstepcounter{Abschnitt}\protect
   \label{#1}\vspace{1ex}
   {\bf (\protect\arabic{section}.\protect\arabic{Abschnitt})}
                     $\qquad$\kkk{#1}}
\newcommand{\zitat}[2]{(\protect\ref{#1}.\protect\ref{#1-#2})}

\begin{document}
\title{The polyhedral Hodge number $h^{2,1}$\\ and vanishing of obstructions}
\author{Klaus Altmann \quad\qquad Duco van Straten }
\date{}
\maketitle
\begin{abstract}
We prove a vanishing theorem for the Hodge number $h^{2,1}$ of 
projective toric varieties provided by a  certain class of polytopes. 
We explain how this Hodge number also gives
information about the deformation theory of the toric Gorenstein
singularity derived from the same polytope. In particular,
the vanishing theorem for $h^{2,1}$ implies that these deformations are
{\em unobstructed}.
\end{abstract}

%
%
\section{Introduction}\label{intro}


\neu{intro-setup} 
For an arbitrary polytope $\kP$, Brion has introduced in \cite{Br} certain
invariants $h^{\ki,\kj}(\kP)$. 
These are defined as dimensions of cohomology groups $H^{\ki,\kj}$ of
complexes which are associated directly to the polytope $\kP$. In case that $\kP$ is a rational polytope, these invariants are exactly the Hodge numbers 
$\dim H^\ki(X,\Omega_{X}^\kj)$
of the corresponding projective toric variety $X=\PP(\kP)$.
\par

In this paper, we focus on the $\kk$-vector spaces $\kT^\kn(\kP):=H^{\kn,1}$    
($\kn \geq 2$) and $\kT^1(\kP):=H^{1,1}/\kk$. The notation is suggested by a
second interpretation of these vector spaces: In \S \ref{deform} we will show that
there is a close relation between $\kT^\kn(\kP)$ and the vector
spaces $T^\kn$ describing the deformation theory of the toric Gorenstein
singularity $X_{\mbox{\rm\kf cone}(\kP)}$ associated to the lattice 
polytope $\kP$.\\

Our main result is a vanishing theorem for $\kT^2(\kP)$ for a certain 
class polytopes. An important special case is
\par

{\bf Theorem} (cf.\ \zitat{van}{pyramid})
{\em
Let $\kP$ be an $n$-dimensional, compact, convex polytope such that every
three-dimensional face is a pyramid.
If no vertex is contained in more than $(n-3)$ two-dimensional,
non-triangular faces, then $\kT^2(\kP)=0$.
}
\par

There is a natural class of polytopes that arises from {\em quivers} (see \cite{AH})
to which these seemingly strange conditions apply. Special examples of such quiver polytopes
appeared in \cite{BCKvS} as a description of toric degenerations of Grassmannians
and partial flag manifolds that appeared in the works of Strumfels \cite{S} and
Lakshmibai \cite{L}. In a forthcoming paper \cite{AvS}, we will apply the above vanishing result to show that the Gorenstein singularities provided by so-called flag like quivers are
{\em unobstructed} and {\em smoothable} in codimension three.   
\par


\neu{intro-contents}
The paper is organized as follows:
\vspace{1ex}\\
In \S \ref{hodge} we recall some notions of homological and
cohomological systems on polyhedral complexes. For the special case
of simplicial sets, these can be found in \cite{GM} or \cite{EMS}. We quote the definition of the 
polyhedral Hodge numbers and review their basic properties.
\vspace{1ex}\\
In \S \ref{Tinv}, we introduce the $\kT$-invariants from a slightly different point of view
as above and show their relation to the polyhedral Hodge numbers. We present some examples
as well as elementary
properties, such as the relation of $\kT^1(\kP)$ to the Minkowski 
decomposition of polytopes.
\vspace{1ex}\\
The paragraphs \S \ref{van} and \S \ref{proof} contain the vanishing theorem for $\kT^2$
and its proof. The result is obtained from a spectral sequence
relating the $\kT$-invariants of a polytope to those of its faces;
$\kT^2(\kP)$ is represented as the kernel of some differential on the
$E_2$-level. In Theorem \zitat{van}{cumbersome}, this description  is transformed
into an explicit set of equations describing $\kT^2(\kP)$.
\vspace{1ex}\\
The final \S \ref{deform} deals with the relations of
the $\kT$-invariants to deformation theory that was mentioned before. 
In the paper \cite{AQ}, a combinatorial description of the cotangent cohomology modules $T^\kn(X_{\mbox{\rm\kf cone}(\kP)})$ was given. From this description it appears that the
$T^{\kn}$ are very sensitive to the interaction of the polytope $\kP$ with the lattice structure of the ambient space. As a consequence, these invariants are often very difficult to calculate explicitly.\\
On the other hand, the invariants $\kT^\kn(\kP)$ are rather coarse; 
they only
depend on the polytope $\kP$ {\em up to projective equivalence}, and the lattice structure is
not involved at all. Nevertheless, in Theorem \zitat{deform}{geq2} we formulate a
sufficient conditions on $\kP$ ensuring that the $T^\kn$ {\em are} determined by
$\kT^\kn(\kP)$. 
In particular, the vanishing Theorem \zitat{van}{pyramid} yields
a vanishing theorem for $T^2(X_{\mbox{\rm\kf cone}(\kP)})$ as well.  
\par


\neu{intro-ack}
{\em Acknowledgement}: We would like to thank M.~Brion and P.~McMullen for
valuable comments and discussions.
\par

%
%
\section{Hodge numbers for polytopes}\label{hodge}


\neu{hodge-cohom}
Let $\cx=\cup_{\kn\geq0}\cx_\kn$ be a finite, {\em polyhedral complex} in a 
$\kk$-vector space $\cxV$ ($\Q\subseteq \kk \subseteq\R$), 
i.e.\ a set of polyhedra in $\cxV$ that is closed under the face operation
and with the additional property 
that for any two $\cxa,\cxb\in\cx$, the intersection
$\cxa\cap\cxb$ is either empty or a common face of both polyhedra.
Here $\cx_\kn$ denotes the subset of $\kn$-dimensional elements of $\cx$.
Examples for such polyhedral complexes 
are {\em simplicial sets} as well as {\em fans}.
\par

{\bf Definition:}
A {\em cohomological system} $\cxF$ on $\cx$ is a covariant functor
from $\cx$ to the category $\AB$ of abelian groups (or to any other
abelian category $\CA$).
\par

Here $\cx$ becomes a small category by declaring
the face relations ``$\cxb\leq\cxa$'' to be the morphisms. So a cohomological
system is nothing else than a collection of abelian groups $\cxF(\cxa)$  
for $\cxa\in\cx$ together       
with compatible face maps $\cxF(\cxb)\to \cxF(\cxa)$.\\
Similarly, a {\em homological system} is defined as a contravariant functor
from $\cx$ to $\AB$.
\par

We fix for each polyhedron $\cxa\in\cx$ an orientation. This enables us to 
introduce for each pair $(\cxb,\cxa)$ of elements of $\cx$ a number
$\keps(\cxb,\cxa)$ as follows:
\kitem
\item
If $\cxb$ is a facet (i.e.\ a codimension-one face) of $\cxa$, then we may
compare the original orientation of $\cxb$ with that inherited from
$\cxa$. Depending on the result, we define
$\keps(\cxb,\cxa):=\pm 1$.
\item
If $\cxb$ is {\em not} a facet of $\cxa$, then we simply set
$\keps(\cxb,\cxa):=0$.
\kenditem
\par

Each cohomological system $\cxF$ on $\cx$ gives rise to a
complex $\cxCF{\kb}$ of abelian groups:
\[
\cxCF{\kn}:= \oplus_{\cxa\in\cx_\kn} \,\cxF(\cxa)
= \oplus_{[\cxa\in\cx;\, \dim \cxa=\kn]} \,\cxF(\cxa)\,.
\]
The differential $d:\cxCF{\kn}\to\cxCF{\kn+1}$ is defined in the obvious way,
using the $\keps(\cxb,\cxa)$ introduced above. The associated
cohomology is denoted by
\[
\cxHF{\kn}:= H^\kn\big(\cxCF{\kb}\big).
\]
Note that there is an analogous construction for homological systems.
\par


\neu{hodge-spectral}
The cohomology groups of a cohomological system $\cxF$
can sometimes be computed
using certain subcomplexes of $\cx$.
To be more precise, 
let $\cxM{\Mi}\subseteq \cx$ be subcomplexes with $\cup_\Mi \cxM{\Mi}=\cx$. 
The {\em nerve} $\nM$ of this covering is the simplicial set defined as
\[
\nM_p:= \{\Mi_0\leq\dots\leq\Mi_p\,|\; 
\cxM{\Mi_0}\cap\dots\cap \cxM{\Mi_p} \neq \emptyset\}.
\]
We obtain cohomological systems
$\CH^q(\cxF)$ on $\nM$ via
\[
\CH^q(\cxF): (\Mi_0\leq\dots\leq\Mi_p) \,\mapsto \,
H^q\big(\cxM{\Mi_0}\cap\dots\cap \cxM{\Mi_p},\, \cxF\big)\,.
\]
\par

{\bf Proposition:}
{\em
There is a degenerating spectral sequence
$E_2^{p,q}=H^p\big(\nM,\,\CH^q(\cxF)\big) \Rightarrow \cxHF{p+q}$
with differentials $d_r: E_r^{p,q}\to E_r^{p+r,q-r+1}$.
}
\par

{\bf Proof:}
Consider the double complex
\[
C^{p,q}:=\oplus_{\cxa\in [\cxM{\Mi_0}\cap\dots\cap \cxM{\Mi_p}]_q} 
\,\cxF(\cxa) \; \mbox{ with }\;
d_I:C^{p,q}\to C^{p+1,q}\,,\; d_{II}:C^{p,q}\to C^{p,q+1}\,.
\]
The first spectral sequence yields 
$E_2^{p,q}=H^p_I H^q_{II}(C^{\kb,\kb})=H^p\big(\nM,\,\CH^q(\cxF)\big)$;
the other one provides the complex $\cxCF{\kb}$ at the $E_1$-level, i.e.\
$\cxHF{\kb}$ is the cohomology of the total complex.
\hfill$\Box$
\par


\neu{hodge-def}
Now assume that $\cf$ is a {\em fan} in the $\cfD$-dimensional
vector space $\cfV$, i.e.\ its
elements are polyhedral cones with $0$ as their common vertex. Note that the
intersection of cones from $\cf$ is always non-empty.
Another special feature of fans is that they come with an important 
cohomological
system for free: $\cfF(\cfa):=\span_k(\cfa)$. 
From this cohomological system ``$\span$'' one derives various other systems
like $\raisebox{0.7ex}{$\cfV$}\hspace{-0.4em}\big/\hspace{-0.2em}
\raisebox{-0.7ex}{$\span$}$ and its exterior powers.
These give rise
to the so called {\em Hodge spaces} of $\cf$, a notion 
which is due to Brion:
\[
H^{\ki,\kj}(\cf):= H^{\cfD-\ki}\Big(\cf,\,\Lambda^\kj 
\raisebox{0.7ex}{$\cfV$}\hspace{-0.4em}\big/\hspace{-0.2em}
\raisebox{-0.7ex}{$\span$} \Big)^\ast\,.
\] 
For rational fans, Danilov has shown in \S 12 of \cite{Da}
that 
$H^{\ki,\kj}(\cf)$
is $H^\ki(X,\Omega_X^\kj)$ where $X=X_\cf$ denotes the toric variety
induced by $\cf$, and $\Omega_X^\kj$ is the reflexive hull of the K\"ahler
$\kj$-differentials on $X_\cf$.
For general fans $\cf$, Brion has obtained the following vanishing results: 
\par

{\bf Proposition:} (cf.\ \S 1 of \cite{Br})
{\em
\kitem
\item[(i)] 
$H^{\ki,\kj}(\cf)=0$ for $\ki<\kj$.
\item[(ii)]
If $|\cf|:=\cup_{
\cfa\in
\cf}\cfa$ 
is not contained in any hyperplane, then $H^{\cfD,\kj}(\cf)=0$
for $\kj<\cfD$, and $H^{\cfD,\cfD}(\cf)$ is isomorphic to $\kk$.
\item[(iii)]
If $|\cf|=\cfV$, and if $e$ is a positive integer such that 
cones with dimension at most $e$ are simplicial, then
$H^{\ki,\kj}(\cf)=0$ for $\ki-\kj>\cfD-e$.
\item[(iv)]
Assume that $|\cf|=\cfV$ and that any two non-simplicial cones in $\cf$
intersect only at the origin. Then $H^{2,1}(\cf)=0$.
\kenditem
}

Note that the assumption of (iv) implies that any $(\cfD-1)$-dimensional cone
in $\cf$ is simplicial. Hence, by (iii), it follows that
$H^{\ki,1}(\cf)=0$ for $\ki\geq 3$.
\par


\neu{hodge-polytope}
Let $\kP\subseteq\kk^n$ be a compact, convex polytope. It gives rise to the
(inner) {\em normal fan} $\fnP$ in the dual space 
$(\kk^n)^\ast\cong\kk^n$. 
Brion has shown
that the diagonal Hodge spaces $H^{\ki,\ki}(\fnP)$ have then a special
combinatorial meaning; they coincide with the 
spaces of the so-called Minkowski $\ki$-weights of $\kP$.\\
In this paper we will focus on the spaces
$H^{\ki,1}(\fnP)$ which sit close to the boundary of the Hodge diamond.
\par

%
%
\section{The $\kT$-invariants}\label{Tinv}


\neu{Tinv-setup}
Let $\kP\subseteq\kk^n$ be a {\em compact, convex polytope}; the cone over
it, denoted by $\cone(\kP)$, generates a (non-complete)
fan $\fcP$ in $\kk^{n+1}$. This gives rise to the following 
invariants of the polytope $\kP$:
\[
\kT^\kn(\kP):= H^\kn\Big(\fcP,\,
\raisebox{0.7ex}{$\kk^{n+1}$}\hspace{-0.4em}\big/\hspace{-0.2em}
\raisebox{-0.7ex}{$\span$} 
\Big) = H^{\kn+1} \big(\fcP,\,\span\big)\,.
\]
The equality is a result of 
the exactness of the complex $C^\kb(\fcP,\,\kk^{n+1})$
sitting in the middle
of the short exact sequence of cohomological systems
\[
0\to \span \longrightarrow \kk^{n+1} \longrightarrow
\raisebox{0.7ex}{$\kk^{n+1}$}\hspace{-0.4em}\big/\hspace{-0.2em}
\raisebox{-0.7ex}{$\span$}
\to 0\,.
\]
Up to isomorphisms, the vector spaces $\kT^\kn(\kP)$ depend only 
on the {\em projective equivalence class} of the given polytope $\kP$.
However, as examples from McMullen and Smilansky show, they are not
combinatorial invariants of $\kP$.
\par

{\bf 
From now on, we will always assume that $\kP\subseteq\kk^n$ has the full
dimension $n$.}
\par

{\bf Lemma:}
{\em
Denote by $\kP^\veee$ the polytope that is {\em polar} to $\kP$, i.e.\ the face
lattice of $\kP^\veee$ is opposite to that of $\kP$, and the cones
$\cone(\kP^\veee)$ and $\cone(\kP)$ are mutually dual. Then, there is a
perfect pairing
\[
\kT^\kn(\kP^\veee) \,\times\, \kT^{n-\kn}(\kP) \longrightarrow \kk\,.
\vspace{-3ex}
\]
}
\par

{\bf Proof:}
If $\cfa\leq\cone(\kP)$ is an $(n+1-\kn)$-dimensional face, then 
$\big[\cfa^\bot\cap\cone(\kP^\veee)\big]\leq \cone(\kP^\veee)$ is a face of 
dimension $\kn$ with
$\span\big[\cfa^\bot\cap\cone(\kP^\veee)\big]=\cfa^\bot$.
Moreover, all faces of $\cone(\kP^\veee)$ arise in this way. Hence,
\[
\renewcommand{\arraystretch}{1.5}
\begin{array}{rcl}
\kT^\kn(\kP^\veee) 
&=&
H^\kn\Big(\fcPd,\,
\raisebox{0.7ex}{$\kk^{n+1}$}\hspace{-0.4em}\big/\hspace{-0.2em}
\raisebox{-0.7ex}{$\span$} \Big)
\;=\;
H^\kn\Big(\fcPd,\, \big((\kb)^\bot\big)^\ast\Big)\\
&=&
H_\kn\big(\fcPd,\, (\kb)^\bot\big)^\ast
\;=\;
H^{n+1-\kn}\big(\fcP,\span\big)^\ast
\;=\;
\kT^{n-\kn}(\kP)^\ast\,.
\vspace{-3ex}
\end{array}
\]
\hspace*{\fill}$\Box$
\par


\neu{Tinv-ex}
The following remarks are intended to obtain a better feeling for the 
meaning of the invariants $\kT^\kn(\kP)$.
\par

\kitem
\item[(i)]
For $\kP=\emptyset$ we define $\cone(\emptyset):=0$, hence
$\kT^\kn(\emptyset)=0$ for every $\kn\in\Z$.
\item[(ii)]
If $\kP$ is a point, then $\cone(\kP)=\kk_{\geq 0}$. In particular, 
$\kT^0(\mbox{point})=\kk$ is the only non-trivial $\kT$-space.
\item[(iii)]
Let $\dim(\kP)\geq 1$. Then, the defining complex for the $\kT^\kn(\kP)$ 
looks like
\[
\begin{array}{ccccccccccccc}
0 & \longrightarrow & C^0 & \longrightarrow & C^1 & \longrightarrow &
\dots & \longrightarrow & C^n & \longrightarrow &
C^{n+1}& \longrightarrow & 0\\
&& || && || &&&& || && ||&&\\
&& \kk^{n+1} && 
\makebox[3em]{$\oplus_{a\in\kP}
\raisebox{0.5ex}{$\kk^{n+1}$}\hspace{-0.3em}\big/\hspace{-0.3em}
\raisebox{-0.5ex}{$\kk\cdot a$}$}
&&&& 
\makebox[3em]{$\oplus_{f<\kP} 
\raisebox{0.5ex}{$\kk^{n+1}$}\hspace{-0.3em}\big/\hspace{-0.2em}
\raisebox{-0.5ex}{$\span\,f$}$}
&& 0&&
\end{array}
\]
with $a\in\kP$ and $f<\kP$ running through the vertices and facets of $\kP$,
respectively. In particular, the injectivity of $C^0\hookrightarrow C^1$
implies $\kT^0(\kP)=0$ and, by the previous lemma, 
$\kT^n(\kP)=\kT^0(\kP^\veee)^\ast=0$.\\
Hence, $\kT^1(\kP),\dots,\kT^{n-1}(\kP)$ are the only non-trivial 
$\kT$-invariants of a polytope $\kP\subseteq\kk^n$.
\kenditem

Denote by $f_j(\kP)$ the number of $j$-dimensional faces of $\kP$
with $f_{-1}:=1$, i.e.\ the Euler equation says $\sum_{j=-1}^n (-1)^j f_j=0$. 
Then $\dim C^\kn=(n+1-\kn)\cdot f_{\kn-1}$.
\vspace{1ex}\\
{\em $n =2$}:\\
The only non-trivial invariant is $\kT^1$ with
$\;\dim \kT^1(\kP) = -\dim C^0 + \dim C^1 - \dim C^2 = f_0(\kP)-3$.
\vspace{1ex}\\
{\em $n =3$}:\\
$\kT^1$ and $\kT^2$ may be non-trivial with
$\;\dim \kT^2(\kP) - \dim \kT^1(\kP) = \sum_\kn (-1)^\kn \dim C^\kn = 
f_2(\kP)-f_0(\kP)$.
\par


\neu{Tinv-brion}
We would like to compare the $\kT$-invariants with Brion's Hodge
spaces. First, there are the straightforward relations
\[
\kT^\kn(\kP)= H^\kn\big(\fcP,\,
\raisebox{0.5ex}{$\kk^{n+1}$}\hspace{-0.4em}\big/\hspace{-0.2em}
\raisebox{-0.5ex}{$\span$} \big)
= H^{n+1-\kn,1}\big(\fcP\big)^\ast
\]
and
\[
\kT^\kn(\kP) = \kT^{n-\kn}(\kP^\veee)^\ast = H^{\kn+1,1}\big(\fcPd\big)\,.
\]
The $\kT$-invariants have also a direct description in terms of
the normal fan $\fnP$ of $\kP$.
\par

{\bf Proposition:}
{\em
Let $\kP\subseteq\kk^n$ be a compact, convex polytope of dimension $n$
and denote by $\fnP$ its inner normal fan.
Then, there is an exact sequence
\[
0\to \kk \longrightarrow H^{1,1}\big(\fnP\big)
\longrightarrow \kT^1(\kP)\to 0\,.
\]
For the remaining indices $\kn\neq 1$ we have 
$H^{\kn,1}\big(\fnP\big)=\kT^\kn(\kP)$. 
}
\par

{\bf Proof:}
Assume that both $\kP$ and $\kP^\veee$ contain the origin as an interior point.
Then, the projection $\kk^{n+1}\surj\kk^n$ induces an isomorphism of fans
$\pi:\partial\,\cone(\kP^\veee)\to\kk_{\geq 0}\cdot\partial\kP^\veee\cong\fnP$.
Moreover, we obtain the following diagram of cohomological systems: 
\[
\dgARROWLENGTH=0.8em
\begin{diagram}
\node[4]{0}
\arrow{s}
\node{0}
\arrow{s}\\
\node[4]{\kk}
\arrow{e,t}{\sim}
\arrow{s}
\node{\kk}
\arrow{s}\\
\node{\mbox{on the fan}\;\partial\,\cone(\kP^\veee):\hspace{3em}}
\node{0}
\arrow{e}
\node{\span}
\arrow{e}
\arrow{s,r}{\sim}
\node{\kk^{n+1}}
\arrow{e}
\arrow{s}
\node{\raisebox{0.4ex}{$\kk^{n+1}$}\hspace{-0.4em}\big/\hspace{-0.2em}
      \raisebox{-0.7ex}{$\span$}}
\arrow{e}
\arrow{s}
\node{0}\\
\node{\mbox{on the fan}\;\fnP:\hspace{5em}}
\node{0}
\arrow{e}
\node{\span}
\arrow{e}
\node{\kk^{n}}
\arrow{e}
\arrow{s}
\node{\raisebox{0.4ex}{$\kk^{n}$}\hspace{-0.4em}\big/\hspace{-0.2em}
      \raisebox{-0.7ex}{$\span$}}
\arrow{e}
\arrow{s}
\node{0}\\
\node[4]{0}
\node{0}
\end{diagram}
\]
Since $H^{\kn,1}\big(\fnP\big)=H^{n-\kn}\big(\fnP,
      \raisebox{0.4ex}{$\kk^{n}$}\hspace{-0.4em}\big/\hspace{-0.2em}
      \raisebox{-0.7ex}{$\span$}\big)^\ast$
and
\[
\kT^\kn(\kP)=\kT^{n-\kn}(\kP^\veee)^\ast=
H^{n-\kn}\big(\fcPd,
      \raisebox{0.4ex}{$\kk^{n+1}$}\hspace{-0.4em}\big/\hspace{-0.2em}
      \raisebox{-0.7ex}{$\span$}\big)^\ast=
H^{n-\kn}\big(\partial\,\cone(\kP^\veee),
      \raisebox{0.4ex}{$\kk^{n+1}$}\hspace{-0.4em}\big/\hspace{-0.2em}
      \raisebox{-0.7ex}{$\span$}\big)^\ast,
\]
the last column of the above diagram implies the long exact sequence
\[
\dots\to H^{n-\kn+1}(\partial\,\cone(\kP^\veee),\kk)^\ast \to
H^{\kn,1}\big(\fnP\big) \to \kT^\kn(\kP) 
\to H^{n-\kn}(\partial\,\cone(\kP^\veee),\kk)^\ast \to\dots\;.
\]
On the other hand, by comparison with the cohomology groups
$H^{\kb}(\fcPd,\kk)=0$, we obtain that 
$H^{\kb}(\partial\,\cone(\kP^\veee),\kk)$ is also trivial -- with the only
exception
$H^{n}(\partial\,\cone(\kP^\veee),\kk)=\kk$.
\hfill$\Box$
\par



\neu{Tinv-T1}
It is well-known that the vector space $H^{1,1}\big(\fnP\big)$ of 
Minkowski $1$-weights is generated (as an abelian group) by the
semi-group of Minkowski summands of $\kk_{\geq 0}$-multiples of $\kP$.
It is useful to see this fact directly:
\par

$H^{1,1}\big(\fnP\big) = H^{n-1}\big(\fnP,
      \raisebox{0.4ex}{$\kk^{n}$}\hspace{-0.4em}\big/\hspace{-0.2em}
      \raisebox{-0.7ex}{$\span$}\big)^\ast
  = H_{n-1}\big(\fnP,(\kb)^\bot\big)$
equals the kernel 
\[
\ker\Big[\oplus_{[\dim \cfa=n-1]}\cfa^\bot \longrightarrow 
\oplus_{[\dim \cfb=n-2]}\cfb^\bot\Big]
=
\ker\Big[\oplus_{d<\kP}\kk\cdot d \longrightarrow
\oplus_{[f<\kP, \,\dim f=2]} \,\span\, f \Big]
\]
with $d\in\kk^n$ running through the edges of $\kP$. The latter space
encodes Minkowski summands of $\kk_{\geq 0}\cdot\kP$ just by 
keeping track of the dilatation
factors of the $\kP$-edges, cf.\ \cite{versal}, Lemma (2.2).
\vspace{1ex}\\
Note that the trivial Minkowski summand $\kP$ itself induces the element
$\ku{1}$ in $H^{1,1}\big(\fnP\big)\subseteq\oplus_{d<\kP}\kk\cdot d$.
It is exactly this element which is killed in the projection
$H^{1,1}\big(\fnP\big)\surj\kT^1(\kP)$ from the previous proposition.
\par

{\bf Corollary:}
{\em
Polytopes $\kP$ with only triangles as two-dimensional faces have a trivial
$\kT^1(\kP)$.
In particular, for simplicial three-dimensional polytopes, 
the only non-trivial $\kT$-invariant is $\kT^2(\kP)$; it has dimension
$f_0(\kP)-4$.
}
\par

{\bf Proof:}
The first claim is clear. The dimension of $\kT^2(\kP)$ for three-dimensional
polytopes follows from \zitat{Tinv}{ex}, the Euler equation, 
and the fact that
$\,3\,f_2(\kP)=2\,f_1(\kP)$ if $\kP$ is simplicial.
\hfill$\Box$
\par

{\bf Examples:} 1) Since the icosahedron $I$ is simplicial, one obtains
$\kT^1(I)=0$ and $\,\dim \kT^2(I)=8$.
\vspace{0.5ex}\\
2) Consider three-dimensional pyramids $P^m$ and 
double pyramids $\kDP^m$ over an $m$-gon; in both cases we have a trivial
$\kT^1$ focusing again the interest on $\kT^2$. Whereas 
$\kT^2(P^m)$ is also trivial, we do have $\,\dim \kT^2(\kDP^m)=m-2$.
\par


\neu{Tinv-double}
We finish this chapter by an extension of the previous example.
We denote by $\kDP(\kP)\subseteq\kk^{n+1}$ the 
{\em double pyramid} over the polytope $\kP\subseteq\kk^n$.
On the polar
level, this means that $\kDP(\kP)^\veee=\kP^\veee\times I$ with
$I:=[0,1]\subseteq\kk^1$. 
\par

{\bf Proposition:}
{\em
The natural inclusion 
$\kT^1(\kP^\veee)\hookrightarrow\kT^1(\kP^\veee\times I)$ has a
one-dimensional cokernel. For $\kn\geq 2$, there are isomorphisms 
$\kT^\kn(\kP^\veee)\stackrel{\sim}{\longrightarrow}\kT^\kn(\kP^\veee\times I)$.
\vspace{0.5ex}\\
Thus, the $\kT$-invariants of a double pyramid depend on those of the base
via 
\[
\kT^\kn\big(\kDP(\kP)\big)=\kT^{\kn-1}(\kP) 
\hspace{0.5em}\mbox{for}\hspace{0.5em} \kn\neq n
\hspace{1em} \mbox{ and}\hspace{1em}
 \dim\, \kT^n\big(\kDP(\kP)\big)=\dim\,\kT^{n-1}(\kP)+1\,.
\vspace{-3ex}
\]
}
\par

{\bf Proof:}
Just to impress the reader,
we are going to use the language of triangulated categories.
The normal fan $\fn(\kP^\veee\times I)$ can be easily expressed by 
$\fn(\kP^\veee)$; if $N,S\in\kk^{n+1}$ denote the ``poles''
$\pm e^{n+1}$, then
\[
\fn(\kP^\veee\times I)= 
\fn(\kP^\veee) \sqcup \fn^N(\kP^\veee) \sqcup \fn^S(\kP^\veee)
\;\mbox{ with }\;
\fn^{N/S}(\kP^\veee):=\big\{\langle\cfa,N/S\rangle\subseteq\kk^{n+1}
      \,|\; \cfa\in\fn(\kP^\veee)\big\}.
\]
Since the complex
$C^\kb\big(\fn(\kP^\veee\times I),\,
\raisebox{0.5ex}{$\kk^{n+1}$}\hspace{-0.4em}\big/\hspace{-0.2em}
\raisebox{-0.4ex}{$\span$} \big)$
is isomorphic to the shifted mapping cone $C^\kb_{(\pi,\pi)}[-1]$ with
\[
\pi\,:\;
C^\kb\Big(\fn(\kP^\veee),\,
\raisebox{0.6ex}{$\kk^{n+1}$}\hspace{-0.4em}\big/\hspace{-0.2em}
\raisebox{-0.6ex}{$\span$} \Big)
\surj
C^\kb\Big(\fn(\kP^\veee),\,
\raisebox{0.6ex}{$\kk^{n}$}\hspace{-0.4em}\big/\hspace{-0.2em}
\raisebox{-0.6ex}{$\span$} \Big)
\]
and $(\pi,\pi):\bullet\to\bullet\oplus \bullet$, 
we obtain that 
$\,C^\kb\big(\fn(\kP^\veee\times I),\,
\raisebox{0.5ex}{$\kk^{n+1}$}\hspace{-0.4em}\big/\hspace{-0.2em}
\raisebox{-0.4ex}{$\span$} \big)[1]$
and
$\,C^\kb\big(\fn(\kP^\veee),\,\kk\big)[1]$ are on top of the distinguished
triangles over the maps $(\pi,\pi)$ and $\pi$, respectively.
Hence, the octahedral axiom for triangulated categories yields a new
distinguished triangle
\[
\dgARROWLENGTH=0.1em
\begin{diagram}
\node[2]{\makebox[3em]{$C^\kb\Big(\fn(\kP^\veee),\,
          \raisebox{0.6ex}{$\kk^{n}$}\hspace{-0.4em}\big/\hspace{-0.2em}
          \raisebox{-0.6ex}{$\span$} \Big)$}}
\arrow{sw,t}{[1]}\\
\node{\,C^\kb\Big(\fn(\kP^\veee),\,\kk\Big)[1]}
\arrow[2]{e}
\node[2]{\,C^\kb\Big(\fn(\kP^\veee\times I),\,
       \raisebox{0.5ex}{$\kk^{n+1}$}\hspace{-0.4em}\big/\hspace{-0.2em}
       \raisebox{-0.4ex}{$\span$} \Big)[1]}
\arrow{nw}
\end{diagram}
\]
inducing the long exact sequence
\[
\dots\to
H^{\kn,1}\big(\fn(\kP^\veee\times I)\big)^\ast
\to
H^{\kn,1}\big(\fn(\kP^\veee)\big)^\ast
\to
H^{n-\kn+2}\big(\fn(\kP^\veee),\,\kk\big)
\to
H^{\kn-1,1}\big(\fn(\kP^\veee\times I)\big)^\ast
\to\dots\;.
\]
As already mentioned at the end of the proof of Proposition
\zitat{Tinv}{brion}, the spaces 
$H^{n-\kn+2}\big(\fn(\kP^\veee),\,\kk\big)$ vanish unless $\kn=2$.
Thus, it remains to use Proposition \zitat{Tinv}{brion} itself, and to
remark that the injection
$\kT^1(\kP^\veee)\hookrightarrow\kT^1(\kP^\veee\times I)$
cannot be an isomorphism.
\hfill$\Box$
\par

{\bf Example:}
Let $\kP$ be the three-dimensional cuboctahedron obtained by cutting the eight 
corners of a cube.
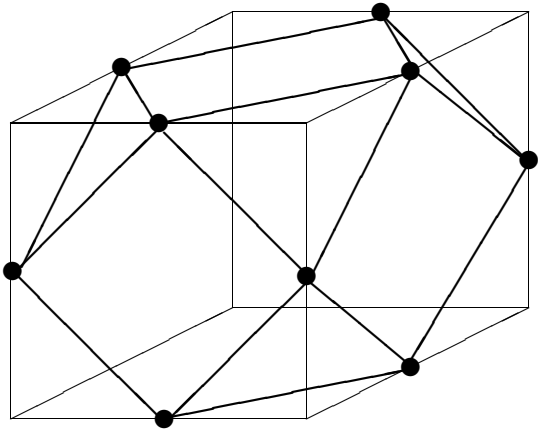
\begin{center}
\unitlength=0.7pt
\begin{picture}(284.00,224.00)(80.00,559.00)
\thinlines
\linethickness{0.1pt}
\put(80.00,560.00){\line(1,0){160.00}}
\put(80.00,720.00){\line(1,0){160.00}}
\put(200.00,780.00){\line(1,0){160.00}}
\put(200.00,620.00){\line(1,0){160.00}}
\put(240.00,560.00){\line(2,1){120.00}}
\put(240.00,720.00){\line(2,1){120.00}}
\put(80.00,720.00){\line(2,1){120.00}}
\put(80.00,560.00){\line(2,1){120.00}}
\put(360.00,620.00){\line(0,1){160.00}}
\put(80.00,560.00){\line(0,1){160.00}}
\put(240.00,560.00){\line(0,1){160.00}}
\put(200.00,620.00){\line(0,1){160.00}}
%
%
\thicklines
\put(161.00,560.00){\line(5,1){136.00}}
\put(297.00,588.00){\line(-6,5){56.00}}
\put(241.00,635.00){\line(-1,-1){76.00}}
\put(360.00,698.00){\line(-3,-5){65.00}}
\put(81.00,640.00){\line(1,-1){81.00}}
\put(82.00,639.00){\line(1,1){78.00}}
\put(160.00,717.00){\line(-3,5){20.00}}
\put(140.00,750.00){\line(-1,-2){54.00}}
\put(161.00,720.00){\line(5,1){137.00}}
\put(298.00,748.00){\line(-1,-2){56.00}}
\put(242.00,636.00){\line(-1,1){79.00}}
\put(298.00,748.00){\line(5,-4){61.00}}
\put(359.00,700.00){\line(-1,1){79.00}}
\put(280.00,779.00){\line(3,-5){18.00}}
\put(142.00,749.00){\line(5,1){139.00}}
\put(163.00,560.00){\circle*{10}}
\put(296.00,588.00){\circle*{10}}
\put(296.00,748.00){\circle*{10}}
\put(140.00,750.00){\circle*{10}}
\put(160.00,720.00){\circle*{10}}
\put(81.00,640.00){\circle*{10}}
\put(240.00,637.00){\circle*{10}}
\put(360.00,700.00){\circle*{10}}
\put(280.00,780.00){\circle*{10}}
\end{picture}
\end{center}
Then, $\kP$ may be decomposed into a Minkowski sum of two tetrahedra --
the summands are formed by the triangles in every other corner. In
particular, $\dim\,\kT^1(\kP)=1$. Moreover, since
$f_0(\kP)=12$, $f_1(\kP)=24$, $f_2(\kP)=14$, 
Example \zitat{Tinv}{ex} implies $\dim\,\kT^2(\kP)=3$.\\
The four-dimensional double pyramid $\kDP(\kP)$ has two kinds of facets:
16 tetrahedra and 12 pyramids over squares.
From the previous proposition, we obtain
$\kT^1\big(\kDP(\kP)\big)=0$, 
$\dim\,\kT^2\big(\kDP(\kP)\big)=1$, and
$\dim\,\kT^3\big(\kDP(\kP)\big)=4$.
\par

%
%
\section{The vanishing theorem}\label{van}


\neu{van-brion}
Let $\kP\subseteq\kk^n$ be an $n$-dimensional, compact, convex 
polytope with $n\geq 1$.
We would like to find conditions under which some of the spaces $\kT^\kn(\kP)$
vanish.\\  
The first idea is to check Brion's properties (iii) and (iv) of
\zitat{hodge}{def} for this purpose. However, we do not
find surprising results in this way -- for instance the first two claims of the
following proposition are already contained in \cite{Br}. The third assertion
generalizes the observation made in Corollary \zitat{Tinv}{T1}.
\par

\begin{minipage}{\textwidth}
{\bf Proposition:}
{\em
\kitem
\item[(1)]
If $\kP$ is a simple polytope, then $\kT^\kn(\kP)=0$ for $\kn\geq 2$. In
particular, for the space of Minkowski summands we obtain
$\,\dim \kT^1(\kP)=\sum_j (-1)^{j+1}\,j\cdot f_j(\kP)$.
\kenditem
}
\end{minipage}
{\em
\kitem
\item[(2)]
If each face of $\kP$ contains at most one non-simple vertex, then we have
still $\kT^2(\kP)=0$.
\item[(3)]
If any $\kl$-face of $\kP$ is a simplex, then $\kT^\kn(\kP)=0$ for $\kn<\kl$.
\kenditem
}
(A vertex of $\kP$ is called {\em simple} if it sits in exactly $n=\dim\kP$
different facets. Moreover, the whole polytope is said to be
simple if every vertex is.)
\par

{\bf Proof:}
The first two claims exploit the fact that 
$\kT^\kn(\kP)=H^{\kn,1}\big(\fnP\big)$
for $\kn \geq 2$. The dimension of $\kT^1(\kP)$ follows from
$\,\dim \kT^1(\kP)=\sum_\kn (-1)^{\kn+1}\dim \kT^\kn (\kP)=
\sum_\kn (-1)^{\kn+1}(n+1-\kn)\cdot f_{\kn-1}$ 
as in \zitat{Tinv}{ex}.\\ 
On the other hand, the proof of the third assertion uses 
\zitat{hodge}{def}(iii) for the ``dual'' fan $\fn(\kP^\veee)$: \\
$\kT^\kn(\kP)=\kT^{n-\kn}(\kP^\veee)^\ast=
H^{n-\kn,1}\big(\fn(\kP^\veee)\big)^\ast=0\,$ if $\,(n-\kn)-1>n-(\kl+1)$.
\hfill$\Box$
\par


\neu{van-spectral}
More interesting results can be obtained from working with the spectral sequence
introduced in \zitat{hodge}{spectral}. When applied to the ``affine'' fan
$\cf:=\fcP$, it provides us with  a nice tool for studying the spaces $\kT^\kn(\kP)$
up to a certain bound $\kn<\kl$:
\par

Write $\kPM{\Mi}<\kP$ for
the $\kl$-dimensional faces and denote
by $\kP^{(\kl)}:=\cup_\Mi\kPM{\Mi}$ the $\kl$-skeleton of our polytope $\kP$.
Then, the simplicial complex $\nM$ with
$\nM_p:= \{\Mi_0\leq\dots\leq\Mi_p\}$
carries the cohomological system
\[
\kTsh^q: (\Mi_0\leq\dots\leq\Mi_p) \,\mapsto \,
\kT^q\big(\kPM{\Mi_0}\cap\dots\cap \kPM{\Mi_p}\big)\,.
\]
\par

{\bf Proposition:}
{\em
There is a degenerating spectral sequence 
with differentials $d_r: E_r^{p,q}\to E_r^{p+r,q-r+1}$
such that
$E_2^{p,q}=H^p\big(\nM,\,\kTsh^q\big) \Rightarrow \kT^{p+q}(\kP)$ for
$p+q<\kl$.
}
\par

{\bf Proof:}
Let $\cf:=\fcP^{(\kl+1)}$
be the union of cones with dimension at most $(\kl+1)$. Then,
using the cohomological system
\[
\CH^q(\span\,\cone):
(\Mi_0\leq\dots\leq\Mi_p) \,\mapsto \, 
H^q\big(\cone\,\kPM{\Mi_0}\cap\dots\cap \cone\,\kPM{\Mi_p},\,\span\big)\,,
\]
section \zitat{hodge}{spectral} yields a degenerating spectral sequence with
\[
E_2^{p,q}=H^p\big(\nM,\,\CH^q(\span\,\cone)\big) \Rightarrow 
H^{p+q}\big(\cf,\,\span\big).
\]
On the other hand, by \zitat{Tinv}{setup} we have
$H^q\big(\cone\,\kPM{\Mi_0}\cap\dots\cap \cone\,\kPM{\Mi_p},\,\span\big)=
\kT^{q-1}\big(\kPM{\Mi_0}\cap\dots\cap \kPM{\Mi_p}\big)$ and,
for $p+q<\kl+1$,
$H^{p+q}\big(\cf,\,\span\big)=H^{p+q}\big(\fcP,\,\span\big)=
\kT^{p+q-1}(\kP)$. Hence, an index shift by one completes the proof.
\hfill$\Box$
\par


\neu{van-supp}
The cohomology groups $E_2^{p,q}=H^p\big(\nM,\,\kTsh^q\big)$ remain unchanged
when we build the complex $C^\kb\big(\nM,\,\kTsh^q\big)$ only from 
the strict tuples $(\Mi_0<\dots<\Mi_p)$. 
\vspace{0.5ex}
In particular, besides $E_2^{0,0}=0$, we obtain at the first glance that 
$E_2^{p,q}=0$ for $q\geq \kl$ or $p\geq 1\,,q=\kl-1$.
However, since the previous proposition restricts us to the region
$p+q<\ell$ anyway, these first vanishings do not help.
\ktrash{\begin{center}
{
\unitlength=1.200000pt
\begin{picture}(230.00,80.00)(-60.00,0.00)
\put(148.00,68.00){\makebox(0.00,0.00){$d_2$}}
\put(128.00,68.00){\vector(2,-1){40.00}}
\put(-60.00,38.00){\makebox(0.00,0.00){
  $\begin{array}{c}E_2^{\kb,\kb}\\(\mbox{for}\, \kl=3)\end{array}$}}
\put(100.00,0){\makebox(0.00,0.00){$p$}}
\put(3.00,80.00){\makebox(0.00,0.00){$q$}}
\put(-5.00,8.00){\makebox(0.00,0.00){$\kT^0$}}
\put(-5.00,28.00){\makebox(0.00,0.00){$\kT^1$}}
\put(-5.00,48.00){\makebox(0.00,0.00){$\kT^2$}}
\put(-5.00,68.00){\makebox(0.00,0.00){$\kT^3$}}
\multiput(2,64)(14,-14){5}{\line(1,-1){9.00}}
\put(8.00,8.00){\circle*{2.00}}
\put(28.00,8.00){\circle*{2.00}}
\put(48.00,8.00){\circle*{2.00}}
\put(68.00,8.00){\circle*{2.00}}
\put(88.00,8.00){\circle*{2.00}}
\put(8.00,28.00){\circle*{2.00}}
\put(28.00,28.00){\circle*{2.00}}
\put(48.00,28.00){\circle*{2.00}}
\put(68.00,28.00){\circle*{2.00}}
\put(88.00,28.00){\circle*{2.00}}
\put(8.00,48.00){\circle*{2.00}}
\put(28.00,48.00){\circle*{2.00}}
\put(48.00,48.00){\circle*{2.00}}
\put(68.00,48.00){\circle*{2.00}}
\put(88.00,48.00){\circle*{2.00}}
\put(8.00,68.00){\circle*{2.00}}
\put(28.00,68.00){\circle*{2.00}}
\put(48.00,68.00){\circle*{2.00}}
\put(68.00,68.00){\circle*{2.00}}
\put(88.00,68.00){\circle*{2.00}}
\put(8.00,8.00){\line(1,0){90.00}}
\put(8.00,78.00){\line(0,-1){70.00}}
\put(8.00,8.00){\circle{5.00}}
\put(28.00,48.00){\circle{5.00}}
\put(48.00,48.00){\circle{5.00}}
\put(68.00,48.00){\circle{5.00}}
\put(88.00,48.00){\circle{5.00}}
\put(8.00,68.00){\circle{5.00}}
\put(28.00,68.00){\circle{5.00}}
\put(48.00,68.00){\circle{5.00}}
\put(68.00,68.00){\circle{5.00}}
\put(88.00,68.00){\circle{5.00}}
\end{picture}}
\end{center}
}
\par

For vertices $a\in\kP$ we denote by $\kP(a)$ the corresponding vertex figure;
it is the polytope obtained by cutting $\kP$ with a hyperplane
sufficiently close to $a$. The faces of $\kP(a)$ are exactly the vertex
figures of those $\kP$-faces containing $a$.
\par

{\bf Lemma:}
{\em
Unless $p=\kl$, the vector spaces $E_2^{p,0}$ on the buttom row vanish.
The remaining one may be expressed by singular cohomology groups with
values in $\kk$ as
\[
E_2^{\kl,0}=\oplus_{a\in\kP\,\mbox{\rm \small vertex}}
\,H^\kl\big(\kP(a), \,\kP(a)^{(\kl-1)}\big) =
\oplus_{a\in\kP\,\mbox{\rm \small vertex}}
\,\widetilde{H}^{\kl-1}\big(\kP(a)^{(\kl-1)}\big)
\]
with $\kP(a)^{(\kl-1)}$ denoting the $(\kl-1)$-skeleton of the 
vertex figure $\kP(a)$.
}
\par

{\bf Proof:}
According to the remark at the beginning of the present section
\zitat{van}{supp}, the vector space
$E_2^{p,0}=H^p\big(\nM,\kTsh^0\big)$ is the $p$-th cohomology of the complex
\[
0\longrightarrow
\oplus_{\Mi_0}\kT^0(\kPM{\Mi_0})
\longrightarrow
\oplus_{\Mi_0<\Mi_1}\kT^0(\kPM{\Mi_0}\cap\kPM{\Mi_1})
\longrightarrow
\oplus_{\Mi_0<\Mi_1<\Mi_2}\kT^0(\kPM{\Mi_0}\cap\kPM{\Mi_1}\cap\kPM{\Mi_2})
\longrightarrow\dots
\]
which will also be denoted by $\kTsh^0$ (creating a slight abuse of
notation).
On the other hand,
for any vertex $a\in\kP$ we call $\kTsh^0(a)$ the complex built similarly as
$\kTsh^0$, but using only those faces $\kPM{\Mi}<\kP$ containing $a$. 
Since $\kT^0$ is trivial unless its argument is a point, the canonical 
projection 
\[
\kTsh^0\surj\oplus_{a\in\kP}\kTsh^0(a)
\]
yields an isomorphism of complexes.\\
%
This splitting enables us to
fix an arbitrary vertex $a\in\kP$ and to forget 
about faces $\kPM{\Mi}$ which do not contain $a$. Then, using the 
vertex figures $\kPM{\Mi}(a)<\kP(a)$, the whole story may be translated
into singular cohomology with values in $\kk$ via
\[
\kT^0\big(\kPM{\Mi_0}\cap\dots\cap\kPM{\Mi_p}\big)
=H^0\big(\kP(a),\, \kPM{\Mi_0}(a)\cap\dots\cap\kPM{\Mi_p}(a)\big).
\]
Denoting by $C_q(\kb)$ the singular $q$-chains,
the Mayer-Vietoris spectral sequence yields
\[
\begin{array}{r@{\;}c@{\;}l}
H^p\big(\kTsh^0(a)\big)= H^p\Big(\Big[
\raisebox{0.5ex}{$C_\kb\big(\kP(a)\big)$}
\hspace{-0.1em}\Big/\hspace{-0.1em}
\raisebox{-0.7ex}{$\sum_\Mi C_\kb\big(\kPM{\Mi}(a)\big)$}
\Big]^\ast\Big)
&=&H^p\big(\kP(a),\,\cup_\Mi \kPM{\Mi}(a)\big)\\
&=&H^p\big(\kP(a),\,\kP(a)^{(\kl-1)}\big)
= \widetilde{H}^{p-1}\big(\kP(a)^{(\kl-1)}\big),
\end{array}
\]  
cf.\ \zitat{proof}{sing} for more details.
Now, the observation that 
the latter groups vanish unless $p=\kl$ finishes the proof.
\hfill$\Box$
\par

{\bf Corollary:}
{\em
(1) Let $\kl\geq 2$. If any at most $\kl$-dimensional face $\kPM{}<\kP$  
satisfies $\kT^\kn(\kPM{})=0$ for $0<\kn<\kl$, then
so does the polytope $\kP$ itself.
\vspace{0.5ex}\\
(2) If there is an $\kl\geq 2$ such that $\kT^1(\kPM{})=0$ for every
$\kl$-face $\kPM{}<\kP$, then $\kT^1(\kP)=0$.
}
\par

{\bf Proof:}
(1) This generalization of the last claim of 
Proposition \zitat{van}{brion}
follows directly from the spectral sequence \zitat{van}{spectral}: The 
assumption means that $E_2^{p,q}=0$ for $p+q<\kl,\, q\neq 0$, and the 
previous lemma takes care of the case $q=0$.
\vspace{0.5ex}\\
(2) Here, the assumption translates into the vanishing $E_2^{0,1}=0$.
\hfill$\Box$
\par


 
\neu{van-T1T2}
For the rest of this chapter, we focus on the case $\kl=3$, i.e.\ we would like
to investigate $\kT^1(\kP)$ and $\kT^2(\kP)$ by studying the 3-dimensional faces
of $\kP$.
Here comes the actual situation of the second layer of our spectral 
sequence (the big circles stand for the vanishing of the corresponding
$E_2$-term):
\vspace{-2ex}
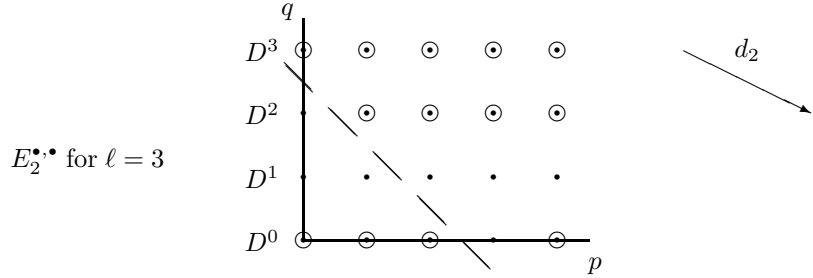
\begin{center}
{
\unitlength=1.200000pt
\begin{picture}(230.00,80.00)(-60.00,0.00)
\put(148.00,68.00){\makebox(0.00,0.00){$d_2$}}
\put(128.00,68.00){\vector(2,-1){40.00}}
\put(-60.00,33.00){\makebox(0.00,0.00){$E_2^{\kb,\kb}\;\mbox{for}\;\kl=3$}}
\put(100.00,0){\makebox(0.00,0.00){$p$}}
\put(3.00,80.00){\makebox(0.00,0.00){$q$}}
\put(-5.00,8.00){\makebox(0.00,0.00){$\kT^0$}}
\put(-5.00,28.00){\makebox(0.00,0.00){$\kT^1$}}
\put(-5.00,48.00){\makebox(0.00,0.00){$\kT^2$}}
\put(-5.00,68.00){\makebox(0.00,0.00){$\kT^3$}}
\multiput(2,64)(14,-14){5}{\line(1,-1){9.00}}
\put(8.00,8.00){\circle*{2.00}}
\put(28.00,8.00){\circle*{2.00}}
\put(48.00,8.00){\circle*{2.00}}
\put(68.00,8.00){\circle*{2.00}}
\put(88.00,8.00){\circle*{2.00}}
\put(8.00,28.00){\circle*{2.00}}
\put(28.00,28.00){\circle*{2.00}}
\put(48.00,28.00){\circle*{2.00}}
\put(68.00,28.00){\circle*{2.00}}
\put(88.00,28.00){\circle*{2.00}}
\put(8.00,48.00){\circle*{2.00}}
\put(28.00,48.00){\circle*{2.00}}
\put(48.00,48.00){\circle*{2.00}}
\put(68.00,48.00){\circle*{2.00}}
\put(88.00,48.00){\circle*{2.00}}
\put(8.00,68.00){\circle*{2.00}}
\put(28.00,68.00){\circle*{2.00}}
\put(48.00,68.00){\circle*{2.00}}
\put(68.00,68.00){\circle*{2.00}}
\put(88.00,68.00){\circle*{2.00}}
\put(8.00,8.00){\line(1,0){90.00}}
\put(8.00,78.00){\line(0,-1){70.00}}
\put(8.00,8.00){\circle{5.00}}
\put(28.00,8.00){\circle{5.00}}
\put(48.00,8.00){\circle{5.00}}
\put(88.00,8.00){\circle{5.00}}
\put(28.00,48.00){\circle{5.00}}
\put(48.00,48.00){\circle{5.00}}
\put(68.00,48.00){\circle{5.00}}
\put(88.00,48.00){\circle{5.00}}
\put(8.00,68.00){\circle{5.00}}
\put(28.00,68.00){\circle{5.00}}
\put(48.00,68.00){\circle{5.00}}
\put(68.00,68.00){\circle{5.00}}
\put(88.00,68.00){\circle{5.00}}
\end{picture}}
\end{center}
\par

{\bf Proposition:}
{\em
Denote by $\kPM{\Mi}<\kP$ 
the three-dimensional faces of $\kP$. Then
\vspace{-1ex}
\kitem
\item[(1)]
$\,\kT^1(\kP)=\ker\Big[ 
\oplus_\Mi \kT^1(\kPM{\Mi})\longrightarrow 
\oplus_{\Mi<\Mj} \kT^1(\kPM{\Mi}\cap\kPM{\Mj})\Big]\,$ and
\vspace{-1ex}
\item[(2)]
if $\,\kT^2(\kPM{\Mi})=0$ for every $\Mi$, then
$\,\kT^2(\kP)=\ker\Big[ d_2:E_2^{1,1}\longrightarrow E_2^{3,0}\Big]\,$. 
\kenditem
}
\par

{\bf Proof:}
The claims follow from
$\kT^1(\kP)=E_{\infty}^{0,1}=E_2^{0,1}=H^0\big(\nM,\kTsh^1\big)$
and, since $E_2^{0,2}=0$ in (2), from
$\kT^2(\kP)=E_{\infty}^{1,1}=E_3^{1,1}$.
\hfill$\Box$
\par


\neu{van-cumbersome}
We are going to apply the previous properties to obtain an explicit
description of
$\kT^2(\kP)$ by equations. In the following we will use the symbols
$a$, $\kPV{}$, $\kPP{}{}$, $\kPF$ to denote $\kP$-faces of dimension 0, 2, 3,
and 4, respectively.
\par

{\bf Notation:}
Whenever $(\kPV{},\kPF)$ is a flag with dimension vector $(2,4)$, then we
denote by $\kPP{(\kPV{},\kPF)}{}$ and $\kPP{}{(\kPV{},\kPF)}$ 
the two unique three-dimensional faces
sitting in between. Their order depends on the orientation of the whole
configuration.\\
For any two-dimensional face $\kPV{}\leq\kP$ we fix some three-dimensional
face $\kPP{}{}(\kPV{})$ containing $\kPV{}$.
\par

{\bf Theorem:}
{\em
Assume that $\kT^1(\kPP{}{})=\kT^2(\kPP{}{})=0$ for every three-dimensional
face $\kPP{}{}\leq\kP$. Then,
$\kT^2(\kP)\subseteq \kk^{\#\{(0,2,3)-\mbox{\rm\kf flags}\}}$ is given by the
following equations in the variables called
$\kts{a}{\kPP{}{}}{\kPV{}}$:
\kitem
\item[(1)]
If $(a,\kPF)$ is a flag with dimension $(0,4)$, then
\[
\sum_{a\in\kPV{}\subseteq\kPF}
\Big[\kts{a}{\kPP{(\kPV{},\kPF)}{}}{\kPV{}} -
\kts{a}{\kPP{}{(\kPV{},\kPF)}}{\kPV{}}\Big]
=0\,.
\makebox[1em][l]{$\hspace*{4em} (1)_{(a,\kPF)}$}
\]
\item[(2)]
For every flag $(\kPV{},\kPP{}{})$, the coordinates 
$\kts{\kb}{\kPP{}{}}{\kPV{}}$ provide an affine relation among the
vertices of $\kPV{}$, i.e.\
\[
\sum_{a\in\kPV{}} \kts{a}{\kPP{}{}}{\kPV{}}\cdot [a,1]=0\,.
\makebox[1em][l]{$\hspace*{8.6em} (2)_{(\kPV{},\kPP{}{})}$}
\]
\item[(3)]
Finally, for each $(0,2)$-flag $(a,\kPV{})$, we simply have
\[
\kts{a}{\kPP{}{}(\kPV{})}{\kPV{}} =0\,.
\makebox[1em][l]{$\hspace*{10.2em} (3)_{(a,\kPV{})}$}
\vspace{-3ex}
\]
\kenditem
}
\par

Note that the equations $(2)_{(\kPV{},\kPP{}{})}$ imply that we can 
completly forget
about the triangular faces $\kPV{}$; they provide only trivial coordinates
$\kts{\kb}{\kb}{\kPV{}}=0$.
\par

The {\em proof of the previous theorem} consists of a detailed, but
straightforward analysis of the differential map
$d_2:E_2^{1,1}\to E_2^{3,0}$. Since it is quite long and technical,
we postpone these calculations to their own section \S \ref{proof}.
In the rest of \S \ref{van}, we continue with a discussion of the 
consequences and applications.
\par


\neu{van-dim4}
{\bf Corollary:}
{\em
If the polytope $\kP$ is four-dimensional, then
$\kT^2(\kP)\subseteq \kk^{\#\{(0,2)-\mbox{\rm\kf flags}\}}$
is given by the easier equations 
\kitem
\item[(1)]
$\;\sum_{\kPV{}\ni a} s(a,\kPV{})=0\,$
for every vertex $a\in\kP$, and
\item[(2)]
$\;\sum_{a\in\kPV{}} s(a,\kPV{})\cdot [a,1]=0\,$ for the two-dimensional faces
$\kPV{}\leq\kP$.
\kenditem
}
\par

{\bf Proof:}
Since $\kPF=\kP$, we may just set
$\kPP{}{}(\kPV{}):=\kPP{}{(\kPV{},\kP)}$ and
$s(a,\kPV{}):= \kts{a}{\kPP{(\kPV{},\kP)}{}}{\kPV{}}$.
\hfill$\Box$
\par

{\bf Example:}
Consider the double pyramid $\kDP(\kP)$ of Example \zitat{Tinv}{double}.
A non-trivial element of the one-dimensional $\kT^2\big(\kDP(\kP)\big)$ 
may be obtained by assigning $\pm 1$ to the vertices of each rectangle
such that adjacent vertices obtain opposite signs.
\par


\neu{van-pyramid}
The main point of the present paper is to provide a vanishing theorem for
$\kT^2(\kP)$ for polytopes
whose three-dimensional faces are not assumed to be simplices.
\par

{\bf Definition:}
We define an inductive process of {\em``cleaning'' vertices and 
two-dimensional faces} of $\kP$. At the beginning, all faces are 
assumed to be ``contaminated'', but then one may repeatedly apply the 
following rules (i) and (ii) in an arbitrary order:
\kitem
\item[(i)]
A two-dimensional $\an$-gon $\kPV{}<\kP$ is said to be clean if at least
$(\an-3)$ of its vertices are so. (In particular, every triangle is
automatically clean.)
\item[(ii)]
A vertex of $\kP$ is declaired to be clean if it is contained in no more than
$(n-3)$ two-dimensional faces that are not cleaned yet.
\kenditem

{\bf Examples:}
(1)
If no vertex of $\kP$ is contained in more than $(n-3)$ two-dimensional,
non-triangular faces, then every vertex and every two-dimensional face
may be cleaned.
\vspace{1ex}\\
(2)
Each vertex of the four-dimensional double pyramid $\kDP(\kP)$ shown in
Example \zitat{Tinv}{double} sits in exactly $2=(n-3)+1$
quadrangular, two-dimensional faces. In particular, it is not possible
to clean any of them at all.
\par

{\bf Theorem:}
{\em
Let $\kP$ be an $n$-dimensional, compact, convex polytope such that 
every three-dimensional face is a pyramid.
If every vertex 
(or, equivalently, every two-dimensional face)
may be cleaned in the sense of the previous definition, 
then $\kT^2(\kP)=0$.
}
\par

{\bf Remarks:}
(1)
Pyramids are the easiest three-dimensional solids with trivial
$\kT$-invariants. Moreover,
polytopes with only pyramids as three-dimensional faces do naturally arise
from quivers, cf.\ \cite{AvS} for more details.
\vspace{1ex}\\
(2) 
The double pyramid $\kDP(\kP)$ from Example \zitat{Tinv}{double}
has a non-trivial $\kT^2$. This shows that the assumption concerning the cleaning of vertices 
cannot be dropped.
\par

{\bf Proof:}
Using the dictionary
\vspace{-2ex}
\begin{center}
\begin{tabular}{r@{\hspace{1em}$\longleftrightarrow$\hspace{1.2em}}l}
``the vertex $a$ is clean''
&
$\kts{a}{\kPP{}{}}{\kPV{}}=0$ for every $\kPV{},\kPP{}{}$\\
``the 2-face $\kPV{}$ is clean''
&
$\kts{a}{\kPP{}{}}{\kPV{}}=0$ for every $a,\kPP{}{}$,
\end{tabular}
\vspace{-2ex}
\end{center}
the vanishing of $\kT^2(\kP)$ is a consequence of Theorem
\zitat{van}{cumbersome} and the following two facts:
\vspace{1ex}\\
(i) If $\kPV{}$ is an $\an$-gon, then, for any $\kPP{}{}$, the 
equation $(2)_{(\kPV{},\kPP{}{})}$ of Theorem \zitat{van}{cumbersome}
says that the coordinates
$\kts{\kb}{\kPP{}{}}{\kPV{}}$ describe an $(\an-3)$-dimensional vector
space. Hence, if $(\an-3)$ of them vanish, then they do all.
In particular, as already mentioned in \zitat{van}{cumbersome}, we do not
have to care about triangular faces $\kPV{}$.  
\vspace{1ex}\\
(ii) Assume that $\kPP{A}{}, \kPP{B}{}$ are two pyramids
with common facet $\kPV{}<\kP$.\\
We denote by $\kP(\kPV{})$ the $(n-3)$-dimensional vertex figure
of a slice of $\kP$ transversal to $\kPV{}$.
In particular, the faces of $\kP(\kPV{})$ correspond to those of $\kP$
containing $\kPV{}$:
While $\bar{\kPV{}}:=\kPV{}(\kPV{})=\emptyset$, the two pyramids
turn into vertices
$\kPPb{A}{}:=\kPP{A}{}(\kPV{})$ and
$\kPPb{B}{}:=\kPP{B}{}(\kPV{})$. Moreover, any four-dimensional face
$\kPF<\kP$ containing $\kPV{}$ corresponds to an edge $\kPFb$ in
$\kP(\kPV{})$.\\
The important feature about pyramids as three-dimensional faces is 
the following:
Any two non-triangular, two-dimensional faces of $\kP$ span an at least 
four-dimensional space.
Hence, for any two-dimensional $\kPV{\prime}$, different from
$\kPV{}$, there is {\em at most one}
four-dimensional $\kPF^\prime<\kP$ containing both
$\kPV{}$ and $\kPV{\prime}$.\\
Thus, if there are given $(n-4)$ 
(contaminated) faces
$\kPV{\Vi}$ additional to $\kPV{}$, then they induce at most $(n-4)$
four-dimensional faces $\kPF^\Vi$
in this way.
Since $\dim\,\kP(\kPV{})=n-3$, this means that it is possible
to find a path along the edges of $\kP(\kPV{})$ connecting the vertices
$\kPPb{A}{}$ and $\kPPb{B}{}$, but avoiding $\bar{\kPF}^\Vi$
($\Vi=1,\dots,n-4$).\\
Let us, w.l.o.g., assume that $\kPPb{A}{}$ and
$\kPPb{B}{}$ are directly connected via an edge $\bar{\kPF}$ with
$\kPF$ not containing the $(n-4)$ faces $\kPV{\Vi}\neq\kPV{}$. Hence,
$\kPPb{A}{}=\kPP{(\kPV{},\kPF)}{}$,
$\kPPb{B}{}=\kPP{}{(\kPV{},\kPF)}$, and
in the equation $(1)_{(a,\kPF)}$ of Theorem \zitat{van}{cumbersome}
\[ 
\sum_{a\in\kb\subseteq\kPF}
\Big[\kts{a}{\kPP{(\kb,\kPF)}{}}{\kb} -
\kts{a}{\kPP{}{(\kb,\kPF)}}{\kb}\Big]
=0\,,
\]
we automatically
sum only over $\kPV{}$ itself and, additionally, over two-dimensional faces 
which are already clean.
\hspace*{\fill}$\Box$
\par

%
%
\section{The proof of the $\kT^2$-equations}\label{proof}

Here, we present the proof of Theorem \zitat{van}{cumbersome}.
It consists of a detailed, but
straightforward analysis of the differential map
$d_2:E_2^{1,1}\to E_2^{3,0}$.
\par


\neu{proof-e11}
{\em Describing $E_2^{1,1}$:}\\ 
According to the remark at the beginning of section \zitat{van}{supp}, 
the vector space
$E_2^{1,1}=H^1\big(\nM,\kTsh^1\big)$ equals the kernel
\[
E_2^{1,1}=\ker\Big[
\oplus_{\Mi_0<\Mi_1}\kT^1(\kPM{\Mi_0}\cap\kPM{\Mi_1})
\longrightarrow \oplus_{\Mi_0<\Mi_1<\Mi_2}
\kT^1(\kPM{\Mi_0}\cap\kPM{\Mi_1}\cap\kPM{\Mi_2})\Big]\,.
\]
Denote by $\kPV{1},\dots,\kPV{\Vn}$ the two-dimensional faces of $\kP$ which
are no triangles; 
each $\kPV{\Vi}$ is contained in some three-dimensional faces
$\kPP{0}{\Vi},\dots,\kPP{\Mn_\Vi}{\Vi}$. 
Note that certain $\kPP{}{}$'s might occur in more 
than one of these lists. Nevertheless, 
\[
E_2^{1,1}\;\cong\; 
\oplus_{\Vi=1}^{\Vn} \kT^1\big(\kPV{\Vi}\big)^{\oplus\Mn_\Vi}
\]
with the $\Mi$-th summand in $\kT^1\big(\kPV{\Vi}\big)^{\oplus\Mn_\Vi}$ being
identified with $\kT^1(\kPP{0}{\Vi}\cap\kPP{\Mi}{\Vi})$; the remaining entries
in $\kT^1(\kPP{\Mi_0}{\Vi}\cap\kPP{\Mi_1}{\Vi})$ may be obtained 
in the usual way as differences
from those of $\kT^1(\kPP{0}{\Vi}\cap\kPP{\Mi_1}{\Vi})$ and 
$\kT^1(\kPP{0}{\Vi}\cap\kPP{\Mi_0}{\Vi})$. On the other hand, if the
intersection $\kPM{\Mi_0}\cap\kPM{\Mi_1}$ is less than two-dimensional,
then $\kT^1(\kPM{\Mi_0}\cap\kPM{\Mi_1})=0$, anyway.\\
We choose the special three-dimensional face $\kPP{}{}(\kPV{\Vi})$ 
mentioned in \zitat{van}{cumbersome} to be $\kPP{0}{\Vi}$.
\par


\neu{proof-d2}
{\em Describing $d_2$:}\\
From \zitat{hodge}{spectral} we recall that 
the double complex inducing the spectral sequence we are dealing with, 
looks as follows:
\[
C^{p,q}=\oplus_{A\in [\kPM{\Mi_0}\cap\dots\cap \kPM{\Mi_p}]_q}
\,\span\big(\cone(A)\big) \; \mbox{ with }\;
d_I:C^{p,q}\to C^{p+1,q}\,,\; d_{II}:C^{p,q}\to C^{p,q+1}\,.
\]
We fix one of the two-dimensional faces $\kPV{\Vi}$ and call it
$\kPV{}$. During \zitat{proof}{d2}, we abbreviate the three-dimensional faces
$\kPP{0}{\Vi},\dots,\kPP{\Mn_\Vi}{\Vi}$ containing $\kPV{\Vi}=\kPV{}$
by $\kPP{0}{},\dots,\kPP{\Mn}{}$. The index $\Mi$ will be reserved for
these $\kPP{\Mi}{}$, while $\Mj\notin\{0,\dots,\Mn\}$ 
points to those three-dimensional faces $\kPP{\Mj}{}<\kP$ 
belonging {\em not} to this list.\\
Assume that $\kPV{}$ is an $\an$-gon with vertices $a^\ai$ 
($\ai\in\Z/\an\Z$). Then, by \zitat{Tinv}{T1},
an element of $\kT^1(\kPV{})$ may be represented
as an $\an$-tuple $(t_1,\dots,t_\an)\in\kk^\an$ with $t_\ai$ being the
dilatation factor assigned to the edge $\ko{a^\ai a^{\ai+1}}<\kPV{}$.
In particular, we may start our tour through the double complex with
an
\[
x=(t^1,\dots,t^\Mn)\in\kT^1\big(\kPV{}\big)^{\oplus\Mn}\subseteq E_2^{1,1}
\;\mbox{ with each $t^\Mi$ represented as }\;
t^\Mi=(t^\Mi_1,\dots,t^\Mi_\an)\in\kk^\an.
\]
The corresponding element $x\in C^{1,1}$ looks like
\[
x_{\Mi_0\Mi_1}\big(\ko{a^\ai a^{\ai+1}}\big)=
(t^{\Mi_1}_\ai-t^{\Mi_0}_\ai)\cdot \kov{30}{10}{a^\ai a^{\ai+1}}
\in
\span(a^\ai, a^{\ai+1})
\hspace{1em}\mbox{ with }\;t^0_\ai:=0\,,
\]
and we have to walk through $C^{\kb,\kb}$ along the following path:
\begin{center}
{
\unitlength=1.200000pt
\begin{picture}(200.00,60.00)(-40.00,0.00)
\put(-65.00,40.00){\makebox(0.00,0.00){Differential $d_2$:}}
\put(160.00,30.00){\makebox(0.00,0.00){$d_{II}$}}
\put(160.00,65.00){\makebox(0.00,0.00){$d_{I}$}}
\put(150.00,20.00){\vector(0,1){20.00}}
\put(150.00,60.00){\vector(1,0){20.00}}
\put(50.00,20.00){\makebox(0.00,0.00){$\mapsto$}}
\put(40.00,30.00){\makebox(0.00,0.00){$\uparrow$}}
\put(20.00,50.00){\makebox(0.00,0.00){$\uparrow$}}
\put(30.00,40.00){\makebox(0.00,0.00){$\mapsto$}}
\put(92.00,20.00){\makebox(0.00,0.00){$d_2(x)\in C^{3,0}$}}
\put(40.00,10.00){\makebox(0.00,0.00){$y$}}
\put(-2.00,40.00){\makebox(0.00,0.00){$x\in C^{1,1}$}}
\put(20.00,60.00){\circle{5.00}}
\put(20.00,20.00){\circle*{2.50}}
\put(40.00,20.00){\circle*{4.00}}
\put(60.00,20.00){\circle*{4.00}}
\put(20.00,40.00){\circle*{4.00}}
\put(40.00,40.00){\circle*{4.00}}
\put(60.00,40.00){\circle*{2.50}}
\put(20.00,60.00){\circle*{2.50}}
\put(40.00,60.00){\circle*{2.50}}
\put(60.00,60.00){\circle*{2.50}}
\end{picture}}
\vspace{-2ex}
\end{center}
The components of the image $d_{I}(x)\in C^{2,1}$ vanish unless exactly two 
of the three indices belong to faces $\kPP{\Mi}{}$ containing $\kPV{}$.
In this case, we obtain
\[
\renewcommand{\arraystretch}{1.3}
\begin{array}{rcl}
d_{I}(x)_{\Mi_0\Mi_1\Mj_2}\big(\ko{a^\ai a^{\ai+1}}\big)&=&
x_{\Mi_1\Mj_2}\big(\ko{a^\ai a^{\ai+1}}\big)
-x_{\Mi_0\Mj_2}\big(\ko{a^\ai a^{\ai+1}}\big)
+x_{\Mi_0\Mi_1}\big(\ko{a^\ai a^{\ai+1}}\big)\\
&=& x_{\Mi_0\Mi_1}\big(\ko{a^\ai a^{\ai+1}}\big)\\
&=&(t^{\Mi_1}_\ai-t^{\Mi_0}_\ai)\cdot \kov{30}{10}{a^\ai a^{\ai+1}}
\end{array}
\]
if $(\kPP{\Mi_0}{}\cap \kPP{\Mi_1}{})\cap\kPP{\Mj_2}{}=
\kPV{}\cap\kPP{\Mj_2}{}=\ko{a^\ai a^{\ai+1}}$.\\
Now, we lift this result to an element $y\in C^{2,0}$, i.e.\ we solve
the equation $d_{II}(y)=d_{I}(x)$. Obviously, the following $y$ does the job:
\[
\renewcommand{\arraystretch}{1.3}
y_{\Mi_0\Mi_1\Mj_2}(a^\ai):=
\left\{\begin{array}{rl}
(t^{\Mi_1}_{\ai-1}-t^{\Mi_0}_{\ai-1}) \cdot a^\ai &
\mbox{if } \kPV{}\cap\kPP{\Mj_2}{}=\ko{a^{\ai-1} a^\ai}\\
(t^{\Mi_1}_\ai-t^{\Mi_0}_\ai) \cdot a^\ai &
\mbox{if } \kPV{}\cap\kPP{\Mj_2}{}=\ko{a^\ai a^{\ai+1}}
\end{array}\right.
\] 
and $y_{\dots}(\kb):=0$ for any other constellation. Its image 
$d_{I}(y)\in C^{3,0}$ asks for quadrupels $(\Mi_0,\Mi_1,\Mj_2,\Mj_3)$
with still exactly two indices belonging to $\kPV{}$-solids.
Up to antisymmetric permutation of the four indices, we have
\[
\renewcommand{\arraystretch}{1.3}
d_{I}(y)_{\Mi_0\Mi_1\Mj_2\Mj_3}(a^\ai)=
\left\{\begin{array}{ll}
(t^{\Mi_1}_{\ai-1}-t^{\Mi_0}_{\ai-1}-t^{\Mi_1}_\ai+t^{\Mi_0}_\ai) \cdot a^\ai
& \mbox{if } \kPV{}\cap\kPP{\Mj_2}{}=\ko{a^\ai a^{\ai+1}}\,,\;
  \kPV{}\cap\kPP{\Mj_3}{}=\ko{a^{\ai-1}a^\ai}\\
(t^{\Mi_1}_{\ai-1}-t^{\Mi_0}_{\ai-1}) \cdot a^\ai 
& \mbox{if } \kPV{}\cap\kPP{\Mj_2}{}=\{a^\ai\}\,,\;
  \kPV{}\cap\kPP{\Mj_3}{}=\ko{a^{\ai-1}a^\ai}\\
-(t^{\Mi_1}_\ai-t^{\Mi_0}_\ai) \cdot a^\ai 
& \mbox{if } \kPV{}\cap\kPP{\Mj_2}{}=\ko{a^\ai a^{\ai+1}}\,,\;
  \kPV{}\cap\kPP{\Mj_3}{}=\{a^\ai\}
\end{array}\right.
\]
and zero for any other constellation. 
The element $d_{I}(y)$ represents the cohomology class
\[
d_2(x)\in E_2^{3,0}=H^3\big(\kTsh^0\big)=
                    \oplus_{a\in\kP}H^3\big(\kTsh^0(a)\big)\,.
\]
Hence, the only non-trivial components $d_2(x)(a\in\kP)$ occur for
$a=a^\ai\in\kPV{}$, and they look like $d_{I}(y)(a^\ai)$ in 
the formula above.
\par


\neu{proof-sing}
{\em Transfer from $H^3\big(\kTsh^0(a)\big)$
to singular cohomology:}\\
Let $a\in\kP$ be an arbitrary vertex.
As already indicated in the proof of Lemma \zitat{van}{supp}, we have to
use the Mayer-Vietoris spectral sequence to describe the isomorphism
\[
H_2\big(\kP(a)^{(2)}\big)=H_3\big(\kP(a),\kP(a)^{(2)}\big)
\stackrel{\sim}{\longrightarrow}
H_3\big(\kTsh_0(a)\big)
\]
with $\kTsh_0(a)$ meaning 
the complex built by the same recipe as
$\kTsh^0(a)$ in \zitat{van}{supp}, but using homology
\[
\kT_0\big(\kPM{\Mi_0}\cap\dots\cap\kPM{\Mi_p}\big)
:=H_0\big(\kP(a),\, \kPM{\Mi_0}(a)\cap\dots\cap\kPM{\Mi_p}(a)\big)
\]
instead of $\kT^0$.
Denoting by $C_q(\kb)$ the singular $q$-chains and abbreviating the vertex
figures $\kPM{\Mi}(a)<\kP(a)$ simply by $\kPMb{\Mi}<\bar{\kP}$, 
we define
\[
K_{p,q}:=\oplus_{\Mi_0\leq\dots\leq\Mi_p}
\raisebox{0.5ex}{$C_q(\bar{\kP})$}
\hspace{-0.1em}\Big/\hspace{-0.1em}
\raisebox{-0.7ex}{$C_q(\kPMb{\Mi_0}\cap\dots\cap\kPMb{\Mi_p})$}
\hspace{1em}\mbox{with}\hspace{1em}
\begin{array}[t]{rl}
d_{I}:&K_{p,q}\to K_{p-1,q}\\
d_{II}:&K_{p,q}\to K_{p,q-1}\,.
\end{array}
\]
The spectral sequence obtained by taking the vertical homology first 
yields the complex $\kTsh_0(a)$ as $E^1_{\kb,0}$ and zero elsewhere.
The other one, beginning with the horizontal homology, has
$E^1_{0,\kb}$ as the only entries at the first level. They form
the complex 
$\raisebox{0.3ex}{$C_\kb(\bar{\kP})$}
\hspace{-0.1em}\big/\hspace{-0.1em}
\raisebox{-0.3ex}{$\sum_\Mi C_\kb(\bar{\kPM{\Mi}})$}$
which is quasi\-isomorphic to
$\raisebox{0.3ex}{$C_\kb(\bar{\kP})$}
\hspace{-0.1em}\big/\hspace{-0.1em}
\raisebox{-0.3ex}{$C_\kb(\cup_\Mi\kPMb{\Mi})$}$.\\
Hence, the existence of the isomorphism promised above is clear. However,
we have to understand what the isomorphism really does with
$[\kPFb]\in H_3(\bar{\kP},\cup_\Mi\kPMb{\Mi})$. To see this
we chase $[\kPFb]$ along the following diagram:
{\small
\[
\begin{array}{@{}c@{}c@{}c@{}c@{}c@{}c@{}c@{}c@{}}
H_3(\bar{\kP},\cup_\Mi\kPMb{\Mi})\\
|\\
\raisebox{0.3ex}{$C_3(\bar{\kP})$}
  \hspace{-0.1em}\big/\hspace{-0.1em}
  \raisebox{-0.3ex}{$\Sigma_\Mi C_3(\bar{\kPM{\Mi}})$}
&\leftarrow&
\oplus_{\Mi_0}\raisebox{0.3ex}{$C_3(\bar{\kP})$}
  \hspace{-0.1em}\big/\hspace{-0.1em}
  \raisebox{-0.3ex}{$C_3(\kPMb{\Mi_0})$}\\
\downarrow&&\downarrow\\
\dots&\leftarrow&
\oplus_{\Mi_0}\raisebox{0.3ex}{$C_2(\bar{\kP})$}
  \hspace{-0.1em}\big/\hspace{-0.1em}
  \raisebox{-0.3ex}{$C_2(\kPMb{\Mi_0})$}
&\leftarrow&
\oplus_{\Mi_0\leq\Mi_1}\raisebox{0.3ex}{$C_2(\bar{\kP})$}
  \hspace{-0.1em}\big/\hspace{-0.1em}
  \raisebox{-0.3ex}{$C_2(\kPMb{\Mi_0}\cap\kPMb{\Mi_1})$}\\
&&\downarrow&&\downarrow\\
&&\dots&\leftarrow&
\oplus_{\Mi_0\leq\Mi_1}\raisebox{0.3ex}{$C_1(\bar{\kP})$}
  \hspace{-0.1em}\big/\hspace{-0.1em}
  \raisebox{-0.3ex}{$C_1(\kPMb{\Mi_0}\cap\kPMb{\Mi_1})$}
&\leftarrow&
\oplus_{\Mi_0\leq\Mi_1\leq\Mi_2}\raisebox{0.3ex}{$C_1(\bar{\kP})$}
  \hspace{-0.1em}\big/\hspace{-0.1em}
  \raisebox{-0.3ex}{$C_1(\kPMb{\ku{\Mi}})$}\\
&&&&\downarrow&&\downarrow\\
&&&&\dots&\leftarrow&
\oplus_{\Mi_0\leq\Mi_1\leq\Mi_2}\raisebox{0.3ex}{$C_0(\bar{\kP})$}
  \hspace{-0.1em}\big/\hspace{-0.1em}
  \raisebox{-0.3ex}{$C_0(\kPMb{\ku{\Mi}})$}
&\leftarrow\dots\\
&&&&&&\downarrow\\
&&&&&&\dots&\leftarrow\dots
\end{array}
\vspace{-2ex}
\]
}
\par

Fixing an arbitrary index $\Mi_0:=0$, we arrive at the third row with
$[\partial\kPFb]\in \raisebox{0.3ex}{$C_2(\bar{\kP})$}
  \hspace{-0.1em}\big/\hspace{-0.1em}\raisebox{-0.3ex}{$C_2(\kPMb{0})$}$.
Let $\partial\kPFb=\kPMb{1}\cup\dots\cup\kPMb{l}$ 
and assume that the orientation of the $\kPMb{\Mi}$ is inherited from some
orientation of $\kPFb$. Then, a possible lift to the right is
\[
-[\kPMb{\Mi}]\,\in\; \raisebox{0.5ex}{$C_2(\bar{\kP})$}
  \hspace{-0.1em}\Big/\hspace{-0.1em}
  \raisebox{-0.7ex}{$C_2(\kPMb{0}\cap\kPMb{\Mi})$}\,,
\quad\Mi=1,\dots,l\,.  
\]
Applying the vertical boundary operator and lifting again to the right, 
we obtain
\[
[\kPMb{\Mi}\cap\kPMb{\Mj}]\,\in\;
  \raisebox{0.5ex}{$C_1(\bar{\kP})$}
  \hspace{-0.1em}\Big/\hspace{-0.1em}
  \raisebox{-0.7ex}{$C_1(\kPMb{0}\cap\kPMb{\Mi}\cap\kPMb{\Mj})$}
\]
with $(\kPMb{\Mi},\kPMb{\Mj})$ running through the pairs of 
mutually adjacent faces of $\kPFb$ with $\Mi<\Mj$. 
Our convention is that the edges $[\kPMb{\Mi}\cap\kPMb{\Mj}]$
inherit their orientation from the first argument, i.e.\
$[\kPMb{\Mi}\cap\kPMb{\Mj}]=-[\kPMb{\Mj}\cap\kPMb{\Mi}]$.
\vspace{1ex}
\[
\begin{array}{@{}c@{}c@{}c@{}c@{}c@{}}
\dots\leftarrow&
\oplus_{\Mi_0\leq\Mi_1\leq\Mi_2}\raisebox{0.3ex}{$C_1(\bar{\kP})$}
  \hspace{-0.1em}\big/\hspace{-0.1em}
  \raisebox{-0.3ex}{$C_1(\kPMb{\ku{\Mi}})$}\\
&\downarrow\\
\dots\leftarrow&
\oplus_{\Mi_0\leq\Mi_1\leq\Mi_2}\raisebox{0.3ex}{$C_0(\bar{\kP})$}
  \hspace{-0.1em}\big/\hspace{-0.1em}
  \raisebox{-0.3ex}{$C_0(\kPMb{\ku{\Mi}})$}
&\leftarrow&
\oplus_{\Mi_0\leq\Mi_1\leq\Mi_2\leq\Mi_3}\raisebox{0.3ex}{$C_0(\bar{\kP})$}
  \hspace{-0.1em}\big/\hspace{-0.1em}
  \raisebox{-0.3ex}{$C_0(\kPMb{\ku{\Mi}})$}\\
&\downarrow&&\downarrow\\
&\dots&\leftarrow&
\oplus_{\Mi_0\leq\Mi_1\leq\Mi_2\leq\Mi_3} H_0(\bar{\kP},\,\kPMb{\ku{\Mi}})&
-\!\!\!- H_3\big[\dots\to
\oplus_{\ku{\Mi}} H_0(\bar{\kP},\,\kPMb{\ku{\Mi}})\to\dots\big]
\end{array}
\]
It is easy to apply the boundary operator to the edges
$[\kPMb{\Mi}\cap\kPMb{\Mj}]$; but for doing the last lifting to 
$\oplus_{\Mi_0\leq\Mi_1\leq\Mi_2\leq\Mi_3}\raisebox{0.3ex}{$C_0(\bar{\kP})$}
  \hspace{-0.1em}\big/\hspace{-0.1em}
  \raisebox{-0.3ex}{$C_0(\kPMb{\ku{\Mi}})$}$, 
we have to
introduce for each vertex $\kPXb\in\kPFb$ an auxillary function $\varphi_\kPX$.
Its arguments are triples of $\kPFb$-facets containing $\kPXb$, and it is
determined by the following properties:
\kitem
\item[(i)]
$\varphi_\kPX(\kPMb{\Mi},\kPMb{\Mj},\kPMb{\Mk})$ is antisymmetric in its
arguments.
\item[(ii)]
If any two of the arguments 
intersect only in $\{\kPXb\}$,
then $\varphi_\kPX(\kPMb{\Mi},\kPMb{\Mj},\kPMb{\Mk}):=0$.
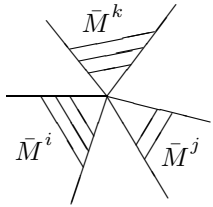
\begin{center}
\unitlength=0.5pt
\begin{picture}(154.00,151.00)(29.00,673.00)
\put(105.00,753.00){\line(-4,5){46.00}}
\put(105.00,753.00){\line(3,-5){46.00}}
\put(105.00,753.00){\line(-1,0){76.00}}
\put(105.00,753.00){\line(-1,-3){27.00}}
\put(105.00,753.00){\line(3,4){52.00}}
\put(105.00,753.00){\line(4,-1){78.00}}
\put(98.00,762.00){\line(5,1){16.00}}
\put(91.00,771.00){\line(5,1){32.00}}
\put(84.00,780.00){\line(5,1){48.00}}
\put(77.00,789.00){\line(5,1){64.00}}
\put(121.00,759.00){\line(-3,-5){7.00}}
\put(132.00,746.00){\line(-3,-5){11.50}}
\put(143.00,743.00){\line(-3,-5){16.00}}
\put(154.00,740.00){\line(-3,-5){20.50}}
\put(88.00,753.00){\line(3,-5){11.00}}
\put(77.00,753.00){\line(3,-5){18.00}}
\put(66.00,753.00){\line(3,-5){25.00}}
\put(55.00,753.00){\line(3,-5){32.00}}
\put(102.00,812.00){\makebox(0,0)[cc]{$\kPMb{\Mk}$}}
\put(161.00,711.00){\makebox(0,0)[cc]{$\kPMb{\Mj}$}}
\put(50.00,715.00){\makebox(0,0)[cc]{$\kPMb{\Mi}$}}
\end{picture}
\end{center}
\item[(iii)]
Denote by $u(\kPXb)$ the number of two-dimensional $\kPFb$-facets meeting in $\kPXb$.
Then, depending on $u(\kPXb)$ and on the fact if there are ``isolated'' 
arguments or not, we distinguish between three cases:
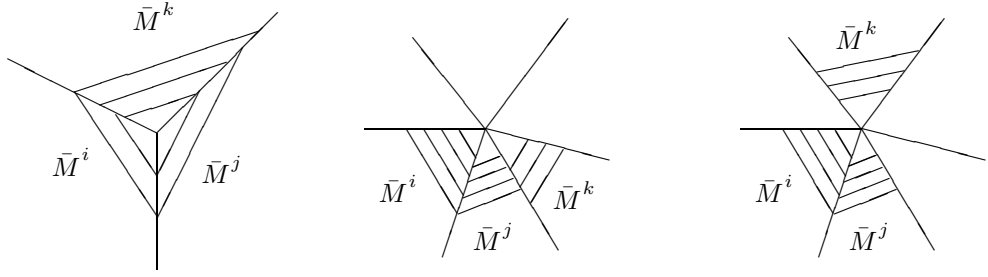
\begin{center}
\hspace*{1em}
\unitlength=0.6pt
\begin{picture}(151.00,203.00)(61.00,620.00)
\put(135.00,750.00){\line(-2,1){94.00}}
\put(135.00,750.00){\line(1,1){76.00}}
\put(135.00,750.00){\line(0,-1){86.00}}
\put(83.00,776.00){\line(3,1){116.00}}
\put(99.00,768.00){\line(3,1){80.00}}
\put(115.00,760.00){\line(3,1){44.00}}
\put(83.00,776.00){\line(2,-3){52.00}}
\put(109.00,762.00){\line(2,-3){26.00}}
\put(162.00,777.00){\line(-1,-2){27.00}}
\put(189.00,804.00){\line(-1,-2){54.00}}
\put(81.00,732.00){\makebox(0,0)[cc]{$\kPMb{\Mi}$}}
\put(176.00,727.00){\makebox(0,0)[cc]{$\kPMb{\Mj}$}}
\put(133.00,824.00){\makebox(0,0)[cc]{$\kPMb{\Mk}$}}
\put(135.00,643.00){\makebox(0,0)[cc]{$u(\kPXb)=3$}}
\put(135.00,615.00){\makebox(0,0)[cc]
    {$\varphi_\kPX(\kPMb{\Mi},\kPMb{\Mj},\kPMb{\Mk}):=1$}}
\end{picture}
\hfill
\begin{picture}(151.00,203.00)(61.00,620.00)
\put(105.00,753.00){\line(-4,5){46.00}}
\put(105.00,753.00){\line(3,-5){46.00}}
\put(105.00,753.00){\line(-1,0){76.00}}
\put(105.00,753.00){\line(-1,-3){27.00}}
\put(105.00,753.00){\line(3,4){52.00}}
\put(105.00,753.00){\line(4,-1){78.00}}
\put(121.00,759.00){\line(-3,-5){7.00}}
\put(132.00,746.00){\line(-3,-5){11.50}}
\put(143.00,743.00){\line(-3,-5){16.00}}
\put(154.00,740.00){\line(-3,-5){20.50}}
\put(88.00,753.00){\line(3,-5){11.00}}
\put(77.00,753.00){\line(3,-5){18.00}}
\put(66.00,753.00){\line(3,-5){25.00}}
\put(55.00,753.00){\line(3,-5){32.00}}
\put(100.00,739.00){\line(5,2){10.00}}
\put(97.00,729.00){\line(5,2){17.00}}
\put(94.00,719.00){\line(5,2){24.00}}
\put(91.00,709.00){\line(5,2){31.50}}
\put(88.00,699.00){\line(5,2){39.00}}
\put(161.00,711.00){\makebox(0,0)[cc]{$\kPMb{\Mk}$}}
\put(50.00,715.00){\makebox(0,0)[cc]{$\kPMb{\Mi}$}}
\put(110.00,685.00){\makebox(0,0)[cc]{$\kPMb{\Mj}$}}
\put(235.00,643.00){\makebox(0,0)[cc]{$u(\kPXb)\geq 4$}}
\put(125.00,615.00){\makebox(0,0)[cc]
    {$\varphi_\kPX(\kPMb{\Mi},\kPMb{\Mj},\kPMb{\Mk}):=2/u(\kPXb)$}}
\end{picture}
\hfill
\begin{picture}(151.00,203.00)(61.00,620.00)
\put(105.00,753.00){\line(-4,5){46.00}}
\put(105.00,753.00){\line(3,-5){46.00}}
\put(105.00,753.00){\line(-1,0){76.00}}
\put(105.00,753.00){\line(-1,-3){27.00}}
\put(105.00,753.00){\line(3,4){52.00}}
\put(105.00,753.00){\line(4,-1){78.00}}
\put(98.00,762.00){\line(5,1){16.00}}
\put(91.00,771.00){\line(5,1){32.00}}
\put(84.00,780.00){\line(5,1){48.00}}
\put(77.00,789.00){\line(5,1){64.00}}
\put(88.00,753.00){\line(3,-5){11.00}}
\put(77.00,753.00){\line(3,-5){18.00}}
\put(66.00,753.00){\line(3,-5){25.00}}
\put(55.00,753.00){\line(3,-5){32.00}}
\put(100.00,739.00){\line(5,2){10.00}}
\put(97.00,729.00){\line(5,2){17.00}}
\put(94.00,719.00){\line(5,2){24.00}}
\put(91.00,709.00){\line(5,2){31.50}}
\put(88.00,699.00){\line(5,2){39.00}}
\put(102.00,812.00){\makebox(0,0)[cc]{$\kPMb{\Mk}$}}
\put(50.00,715.00){\makebox(0,0)[cc]{$\kPMb{\Mi}$}}
\put(110.00,685.00){\makebox(0,0)[cc]{$\kPMb{\Mj}$}}
\put(125.00,615.00){\makebox(0,0)[cc]
    {$\varphi_\kPX(\kPMb{\Mi},\kPMb{\Mj},\kPMb{\Mk}):=1/u(\kPXb)$}}
\end{picture}
\hspace*{1em}
\end{center}
\kenditem
Now, it is not difficult to check that a possible lifting of the tuple
$\big(\partial[\kPMb{\Mi}\cap\kPMb{\Mj}]\big)_{\Mi<\Mj}$
to the vector space
$\oplus_{\Mi_0\leq\dots\leq\Mi_3}\raisebox{0.3ex}{$C_0(\bar{\kP})$}
  \hspace{-0.1em}\big/\hspace{-0.1em}
  \raisebox{-0.3ex}{$C_0(\kPMb{\ku{\Mi}})$}$
is given by
\[
\varphi_\kPX(\kPMb{\Mi},\kPMb{\Mj},\kPMb{\Mk})\cdot [\kPXb]
\;\in\; 
\raisebox{0.7ex}{$C_0(\bar{\kP})$} \hspace{-0.1em}\big/\hspace{-0.1em}
\raisebox{-0.7ex}{$C_0(\kPMb{0}\cap\kPMb{\Mi}\cap\kPMb{\Mj}\cap\kPMb{\Mk})$}
\]
with $(\kPMb{\Mi},\kPMb{\Mj},\kPMb{\Mk})$ running through all triples of 
facets of $\kPFb$ with $\kPMb{\Mi}\cap\kPMb{\Mj}\cap\kPMb{\Mk}=\{\kPXb\}$
and $\Mi<\Mj<\Mk$. The projection to
$H_0(\bar{\kP},\,\kPMb{0}\cap\kPMb{\Mi}\cap\kPMb{\Mj}\cap\kPMb{\Mk})$
does not change the shape of the element 
$\varphi_\kPX(\kPMb{\Mi},\kPMb{\Mj},\kPMb{\Mk})\cdot [\kPXb]$. However, if
$\kPXb\in\kPMb{0}$, then the element $[\kPXb]$ spanning the whole homology 
group vanishes, anyway.\\
Note that the final result of the previous construction cannot depend on
the choice of $\kPMb{0}$ made in the very beginning. 
In particular, one might exploit this freedom to take for $\kPMb{0}$ one
of the $\kPFb$-faces, or to do exactly the opposite.
\par


\neu{proof-interpret}
{\em Interpreting $d_2(x)(a^\ai)$ inside
$H^3\big(\kP(a^\ai),\kP(a^\ai)^{(2)}\big) $:}\\
%
We apply the previous calculation to show how $d_2(x)(a^\ai)$ acts on
a homology class 
$[\kPFb]=[\kPF(a^\ai)]\in H_3\big(\kP(a^\ai),\kP(a^\ai)^{(2)}\big)$ induced
by a four-dimensional face $\kPF<\kP$ containing $a^\ai$:
Unless $\kPV{}\subseteq\kPF$, we have 
$\big\langle d_2(x)(a^\ai), [\kPFb]\big\rangle=0$. However,
assuming $\kPV{}\subseteq\kPF$, then $\kPF$ contains exactly two faces
$\kPP{\Mi}{},\kPP{\Mk}{}$ with common facet $\kPV{}$, and the result is
\[
\big\langle d_2(x)(a^\ai), [\kPFb]\big\rangle \;=\;
  \big(t^{\Mk}_{\ai-1}-t^{\Mk}_\ai\big) - 
  \big(t^{\Mi}_{\ai-1}-t^{\Mi}_\ai\big)\,.
\]
{\em Proof:}
Using the notation
of \zitat{proof}{d2}, we select the first $\kPV{}$-solid 
$\kPP{0}{}=\kPP{}{}(\kPV{})$ as the face inducing the vertex figure
$\kPMb{0}<\bar{\kP}=\kP(a^\ai)$ being fixed in \zitat{proof}{sing}.
Then, the only quadrupels having a chance to produce a non-trivial entry 
in both steps are
\[
[0,\Mi,\Mj_2,\Mj_3] \hspace{1em} \mbox{with }
\renewcommand{\arraystretch}{1.3}
\begin{array}[t]{l@{\hspace{0.5em}}l}
\bullet&
\Mi\in\{0,\dots,\Mn\}
\hspace{0.8em}\mbox{and}\hspace{0.8em}
\kPP{\Mi}{},\,\kPP{\Mj_2}{},\,\kPP{\Mj_3}{}<\kPF \\
\bullet&
\kPV{}\cap\kPP{\Mj_2}{},\, \kPV{}\cap\kPP{\Mj_3}{}=
\{a^\ai\},\; \ko{a^{\ai-1}a^\ai},\; \mbox{ or }\; \ko{a^{\ai}a^{\ai+1}}\\
\bullet&
\kPV{}\cap\kPP{\Mj_2}{}\cap\kPP{\Mj_3}{}=\{a^\ai\}\,.
\end{array}
\]
Focusing on the vertex figures 
$\bar{\kPV{}}\leq\kPPb{\Mi}{}\leq\kPFb$ at $a^\ai$, we see that 
$\kPPb{\Mi}{}$ is a polygon with 
$\bar{\kPV{}}=\ko{{a}^{\ai-1}{a}^{\ai+1}}$ as one of its edges.
While $\bar{\kPV{}}\cap\kPPb{\Mj_2}{}\cap\kPPb{\Mj_3}{}=\emptyset$, the
result of \zitat{proof}{sing} implies that the
intersection $\kPPb{\Mi}{}\cap\kPPb{\Mj_2}{}\cap\kPPb{\Mj_3}{}$ has to be
some point $\kPXb\neq\,\ko{a}^{\ai-1},\,\ko{a}^{\ai+1}$. 
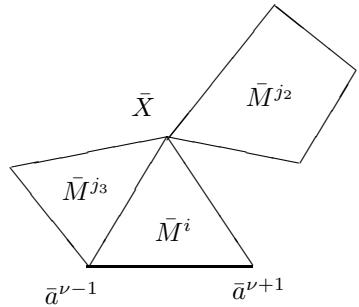
\begin{center}
\unitlength=0.6pt
\begin{picture}(220.00,160.00)(54.00,670.00)
\put(155.00,755.00){\line(-5,-1){98.00}}
\put(155.00,755.00){\line(-3,-5){49.00}}
\put(155.00,755.00){\line(2,-3){54.00}}
\put(155.00,755.00){\line(5,-1){83.00}}
\put(239.00,738.00){\line(3,5){35.00}}
\put(274.00,797.00){\line(-5,4){51.00}}
\put(223.00,838.00){\line(-4,-5){68.00}}
\put(56.00,735.00){\line(4,-5){49.50}}
\put(104.50,673.00){\line(1,0){104.50}}
\put(95.00,656.00){\makebox(0,0)[cc]{$\bar{a}^{\ai-1}$}}
\put(213.00,659.00){\makebox(0,0)[cc]{$\bar{a}^{\ai+1}$}}
\put(140.00,775.00){\makebox(0,0)[cc]{$\kPXb$}}
\put(220.00,779.00){\makebox(0,0)[cc]{$\kPPb{\Mj_2}{}$}}
\put(104.00,720.00){\makebox(0,0)[cc]{$\kPPb{\Mj_3}{}$}}
\put(159.00,694.00){\makebox(0,0)[cc]{$\kPPb{\Mi}{}$}}
\end{picture}
\end{center}
Hence, fixing $\kPPb{\Mi}{}$ and
choosing the ordering of the $\kP$-faces and their corresponding indices 
well, we obtain contributions to 
$\big\langle d_2(x)(a^\ai), [\kPFb]\big\rangle$ only from the arguments
$\kPPb{\Mj_2}{}$ and $\kPPb{\Mj_3}{}$ running through the two-dimensional
$\kPFb$-faces fitting in one of the following cases:
\kitem
\item[(i)]
$\kPPb{\Mj_2}{}$ has a common edge $\ko{a^{\ai+1}\kPX}$ with $\kPPb{\Mi}{}$. 
Then, if $\kPPb{\Mj_3}{}\ni\kPXb$ is adjacent to one of them, we obtain twice
\[
d_2(x)_{0\Mi\Mj_2\Mj_3}(a^\ai)\cdot
\varphi_\kPX(\kPMb{\Mi},\kPMb{\Mj_2},\kPMb{\Mj_3})\;=\;
-t^{\Mi}_\ai \cdot 2/u(\kPXb)\,.
\]
Moreover, there are $\big(u(\kPXb)-4\big)$ possibilities such that 
$\kPPb{\Mj_3}{}\ni\kPXb$ is ``isolated''. Each constellation yields 
the contribution
$-t^{\Mi}_\ai\cdot 1/u(\kPXb)$.
\item[(ii)]
$\kPPb{\Mj_3}{}$ has a common edge $\ko{\kPX a^{\ai-1}}$ with $\kPPb{\Mi}{}$.
Then, as in (i), we obtain twice $t^{\Mi}_{\ai-1}\cdot 2/u(\kPXb)$
and $\big(u(\kPXb)-4\big)$-times $t^{\Mi}_{\ai-1}\cdot 1/u(\kPXb)$.
\item[(iii)]
If both $\ko{a}^{\ai+1}\notin\kPPb{\Mj_2}{}$ and
$\ko{a}^{\ai-1}\notin\kPPb{\Mj_3}{}$, then the result of
\zitat{proof}{d2} shows that this case contributes
nothing to $\big\langle d_2(x)(a^\ai), [\kPFb]\big\rangle$.
\kenditem
Altogether, this adds up to $(t^{\Mi}_{\ai-1}-t^{\Mi}_\ai)$,
and we should finally remark that the exceptional cases ``$u(\kPXb)=3$''
and ``$\kPPb{\Mi}{}$ is a triangle'' yield the same result. In the latter
situation, the cases (i) and (ii) might overlap.
\hfill$(\Box)$
\par

Finally, we should remark that it is exactly the differences
$t^{\Mi}_{\ai-1}-t^{\Mi}_\ai$ which are called
$\kts{a^\ai}{\kPP{\Mi}{}}{\kPV{}}$ in \zitat{van}{cumbersome}.
The equations $(2)_{(\kPV{},\kPP{}{})}$ of the theorem say nothing
else than that these $s$-variables come from some $t$'s satisfying the
equations for Minkowski summands of $\kPV{}$ as mentioned in \zitat{Tinv}{T1}.  
\hspace*{\fill}$\Box$
\par

%
%
\section{Applications to deformation theory}\label{deform}


\neu{deform-setup}
Let $N,M$ be two finitely generated, free abelian groups which are mutually
dual;
denote by $N_\R, M_\R$ the vector spaces obtained by extending the scalars.
Each polyhedral, rational cone $\sigma\subseteq N_\R$ with apex in $0$ 
gives rise to an affine toric variety $X_\sigma:=\Spec\,\C[\sigma^\veee\cap M]$. 
It comes with an action of the torus
$N_\C\otimes_\C\C^\ast=\Spec\, \C[M]$, which leads to a stratification into orbits which are parametrized by the
faces of $\sigma$. We refer to \cite{Da} for more details.\\
In particular, the trivial face $\sigma\leq\sigma$ corresponds to a unique fixed point 
$\mbox{\rm orb}(\sigma)=0$ of
the torus action. It is the ``most singular'' point of $X_\sigma$, and we are
going to study the deformation theory of the germ $(X_\sigma,0)$.
\par

The point that makes toric varieties so exciting is the fact that many
algebro-geometric properties of $X_\sigma$ (or its non-affine
generalizations) translate directly into combinatorial properties of cones
and their relation to the lattice structure $N\subseteq N_\R$. A first
example of such a translation can be seen in \zitat{hodge}{def}. We will need in the future
the following two further examples of such translations:
\kitem
\item
$X_\sigma$ is Gorenstein if and only if 
$\sigma$ is the cone over a compact, convex
lattice polytope $\kP\subseteq\R^n$ sitting in an affine hyperplane of 
height one. 
This means that $N=\Z^n\times\Z$, and
$\kP$ is a polytope with vertices in $\Z^n\times\{1\}$. 
\item
$X_\sigma$ is, additionally, smooth in codimension two iff the edges of $\kP$
do not contain any interior lattice points. 
\kenditem
\par


\neu{deform-Ti}
If $X=\Spec\,A$ is an affine algebraic variety, then the cohomology of the cotangent
complex produces $A$-modules $T^\kn_X$ which play an important role in deformation
theory: $T^0_X$ describes infinitesimal automorphisms, $T^1$ describes infinitesimal
deformations, and $T^2_X$ contains the obstructions for extending
deformations to larger base spaces. 
See \cite{BeCh} for a nice survey, or \cite{Loday} for the details.\\
In the case that $X=X_\sigma$ is a toric variety, the ring
$A=\C\,[\sigma^\veee\cap M]$ as well as the modules $T^\kn_X$ are
$M$-graded. It is possible to obtain combinatorial formulas for the
homogeneous pieces $T^\kn_X(-R)$ with $R\in M$. This has been done in
\cite{AQ}, and we recall the result:
\par

Assume we are given a rational, polyhedral cone 
$\sigma=\langle a^1,\dots,a^\an\rangle \subseteq N_{\R}$
with $a^1,\dots,a^\an\in N$ denoting its {\em primitive
fundamental generators}, 
i.e.\ none of the $a^\ai$ is a proper multiple of an element of $N$.
The dual cone is
$\sigma^{\veee}:= \{ r\in M_{\R}\,|\; \langle \sigma,\,r\rangle \geq 0\} 
\subseteq M_{\R}$. For any degree $R\in M$ and
face $\tau\leq\sigma$ we introduce a special subset of lattice points of $\sigma^{\veee}$:
\[
K_\tau^R:= \sigma^\veee \cap \big(R-\innt\, \tau^\veee\big)\cap M
\subseteq (\sigma^\veee\cap M)\,.
\]
In particular, $K_0^R=\sigma^\veee\cap M$, whereas $K_\sigma^R$ consists of a finite
set of lattice points.
For an arbitrary subset $K \subset M$ we set:
\[
\bHom(K,\C):= \big\{f:K\to\C\,\big|\;
f(r)+f(s)=f(r+s) \;\mbox{ if }\; r,s,r+s\in K\big\}\;.
\]
For each given $R \in M$, these sets give rise to a cohomological system
$\bHom\big(K^R_\kb,\C\big)$ on the ``affine'' fan $\fsigma$.
\par

{\bf Theorem:}
(cf.\ \cite{AQ}, (5.3))\quad
{\em
For  $\kn\leq 2$, one has 
$$T^\kn_X(-R)= H^\kn\big(\fsigma,\,\bHom(K^R_\kb,\C)\big)\;.$$

Moreover, if either $\kn \leq 1$, or if $\kn=2$ and $X_\sigma$ is Gorenstein in
codimension two, then 
$$T^\kn_X(-R)=H^\kn\big(\fsigma,\,\span_\C(K^R_\kb)^\ast\big)\;.$$}
\par

{\bf Remarks:}
1) There is always a natural homomorphism of cohomological systems
$\,\span_\C(K^R_\kb)^\ast \to \bHom(K^R_\kb,\C)$, but in general their 
cohomology groups are different. The second part of the theorem thus gives a
condition under which we can replace the complicated system $\bHom(K^R_\kb,\C)\big)$
by a slightly simpler system of {\em vector spaces}. 
\vspace{1ex}\\
%
2) The module structure of $T^\kn$ is the natural one: If 
$x^s\in\C[\sigma^\veee\cap M]$, then the multiplication with $x^s$ is 
obtained from the map $T^\kn_X(-R)\to T^\kn_X(-R+s)$ provided by the inclusions
$K^{R-s}_\tau\subseteq K^R_\tau$.
\vspace{1ex}\\
3) The property {\em  Gorenstein in codimension two} translates into the following
condition for the cone: For every two-dimensional face 
$\langle a^\ai,a^\aj\rangle<\sigma$ there is an $R_{\ai\aj}\in M$
with $\langle a^\ai,R_{\ai\aj}\rangle = \langle a^\aj,R_{\ai\aj}\rangle =1$.
\par


\neu{deform-leq1}
Let $\kP\subseteq\R^n$ be a lattice polytope; via $\sigma:=\cone(\kP)$
it gives rise to a {\em toric Gorenstein singularity} $X:=X_\kP$. 
For this special case, we are going to 
explain the relations between the vector spaces $T^\kn_X(-R)$ and the 
coarse $\kT$-invariants $\kT^\kn$ defined in \S \ref{Tinv}.
\vspace{1ex}\\
If $a^1,\dots,a^\an\in\Z^n$ denote the vertices of $\kP$, 
then $a^\ai:=(a^\ai,1)\in N$ are the fundamental generators of $\sigma$.
Moreover, there is a special degree $R^\ast:=[\ku{0},1]\in M$; it recovers
the polytope from the cone via $\kP=\sigma\cap [R^\ast=1]$.
\par

{\bf Proposition:} 
{\em
Let $\kP$ and $X:=X_\kP$ be as before. If $R\in M$ is a degree such that
$R\leq 1$ holds everywhere on $\kP$, then
$\kP\cap [R=1]$ is a face of $\kP$ and
\[
T^\kn_X(-R)= \kT^\kn\Big(\kP\cap [R=1]\Big) 
\hspace{1em}\mbox{for}\hspace{0.8em} \kn\leq 2.
\vspace{-3ex}
\]
}
\par

{\bf Proof:}
The reader should convince her/himself from the fact that
the property $R \leq 1$ in $\kP$ imlpies
\[
\span_\C(K^R_\tau)=\left\{\begin{array}{ll}
\tau^\bot & \mbox{if } \tau\leq \cone(\kP\cap [R=1])\\
0 & \mbox{otherwise}\,.
\end{array}\right.
\vspace{-3ex}
\]

The claim then follows from Theorem \zitat{deform}{Ti}. \hspace*{\fill}$\Box$
\par


\neu{deform-loc}
It is possible to describe $T^1_X(-R)$
in the Gorenstein case also for degrees with $R\not\leq 1$ on $\kP$.
However, in the following three sections of the present paper, 
we look for sufficient
conditions forcing $T^1_X(-R)$ and $T^2_X(-R)$ to vanish for those $R$.
\par

Assume that $\sigma=\langle a^1,\dots,a^\an\rangle \subseteq N_\R$
is a rational, polyhedral cone as in \zitat{deform}{Ti}. 
For any degree $R\in M$, we define another
homological system $V^R_\kb \supseteq\span_\C(K^R_\kb)$ on $\fsigma$ by
\[
\renewcommand{\arraystretch}{1.1}
V^R_\tau:=\bigcap_{a^\ai\in\tau} V^R_{a^\ai}
\hspace{2em}\mbox{with}\hspace{1.5em}
V^R_{a^\ai}:=\span_\C(K^R_{a^\ai})=
\;\left\{
\begin{array}{@{}ll}
M_\C & \mbox{if } \,\langle a^\ai,R\rangle\geq 2\\
(a^\ai)^\bot & \mbox{if } \,\langle a^\ai,R\rangle=1\\
0 & \mbox{if } \,\langle a^\ai,R\rangle\leq 0\,.
\end{array}\right.  
\vspace{-2ex}
\]
\par

Let $X_\sigma$ be {\em smooth in codimension two}, i.e.\ 
whenever $\langle a^\ai,a^\aj\rangle<\sigma$ is a two-dimensional face, 
then the set $\{a^\ai,a^\aj\}$ may be extended to a $\Z$-basis of the 
whole lattice $N$. In particular, for any $R\in M$, we have 
$V^R_{\langle a^\ai,a^\aj\rangle} =
\span_\C(K^R_{\langle a^\ai,a^\aj\rangle})$
for these faces. Hence, for $X_{sigma}$ smooth in codimension two one has
$$T^1_X(-R)=H^1\big(\fsigma, \,(V^R_\kb)^\ast\big) .$$
\par

{\bf Definition:}
If $X_\sigma$ is smooth in codimension two, then we define the
{\em local contribution} of
a three-dimensional face $\tau\leq\sigma$ to $T^2_X(-R)$ as
\[
T^2_{\tau,\,\mbox{\kf loc}}(-R)
\;:=\;
\left(
\raisebox{1ex}{
$V^R_\tau$
}\hspace{-0.5em}\Big/\hspace{-0.5em} \raisebox{-1ex}{
$\span_\C K^R_\tau$}\right)^\ast
=
\left(
\raisebox{1ex}{
$\bigcap_{a^\ai\in\tau}\big(\span_\C K^R_{a^\ai}\big)$
}\hspace{-0.5em}\Big/\hspace{-0.5em} \raisebox{-1ex}{
$\span_\C \big(\bigcap_{a^\ai\in\tau}K^R_{a^\ai}\big)$}
\right)^\ast
\,.
\]
\par

If $\dim \sigma=3$ itself, then Theorem \zitat{deform}{Ti} tells us that
$T^2_X(-R)=T^2_{\sigma,\,\mbox{\kf loc}}(-R)$. Moreover, 
for the general case, we obtain the straightforward
\par

{\bf Proposition:}
{\em
Let $X_\sigma$ be smooth in codimension two. If there are no local
contributions from three-dimensional faces to $T^2_X(-R)$
(i.e.\ if $\,T^2_X$ sits in codimension at least four), then
$$T^2_X(-R)=H^2\big(\fsigma,\,(V^R_\kb)^\ast\big)\; .$$
}
\par

{\bf Application:}
If the three-dimensional faces of $\sigma$ are either smooth (generated by a
part of a $\Z$-basis of $N$) or isomorphic to cones over unit squares in
$\Z^2$, then $X_\sigma$ is a {\em conifold} in codimension three, 
i.e.\ it is smooth in codimension two
and has at most $A_1$-singularities in codimension three. In particular, 
for thoses cones,
the assumption of the previous proposition is satisfied for
every multidegree $R\in M$.
\par

{\bf Example:}
To get some familiarity with the sets $K^R_\tau$, we explain the vanishing of the local 
contributions for conifolds on the combinatorial level.. Let $\tau$ be the cone over
a unit square.
Unless $R$ is positive at the four vertices of this square, the space
$\bigcap_{a^\ai\in\tau}\big(\span_\C K^R_{a^\ai}\big)$ vanishes, anyway.
Now, focusing on these four positive
values, there are only the following possibilities:
\begin{center}
\hspace*{\fill}
\unitlength=0.4pt
\begin{picture}(80.00,80.00)(30.00,730.00)
\put(40.00,800.00){\line(1,0){60.00}}
\put(100.00,800.00){\line(0,-1){60.00}}
\put(100.00,740.00){\line(-1,0){60.00}}
\put(40.00,740.00){\line(0,1){60.00}}
\put(40.00,800.00){\circle*{6}}
\put(100.00,800.00){\circle*{6}}
\put(100.00,740.00){\circle*{6}}
\put(40.00,740.00){\circle*{6}}
\put(20.00,820.00){\makebox(0,0)[cc]{${\ks \langle a^4,R\rangle=1}$}}
\put(122.00,820.00){\makebox(0,0)[cc]{${\ks\langle a^3,R\rangle=1}$}}
\put(20.00,720.00){\makebox(0,0)[cc]{${\ks \langle a^1,R\rangle=1}$}}
\put(122.00,720.00){\makebox(0,0)[cc]{${\ks\langle a^2,R\rangle=1}$}}
\put(70.00,770.00){\makebox(0,0)[cc]{$\tau$}}
\end{picture}
\hspace*{\fill}
\unitlength=0.4pt
\begin{picture}(80.00,80.00)(30.00,730.00)
\put(40.00,800.00){\line(1,0){60.00}}
\put(100.00,800.00){\line(0,-1){60.00}}
\put(100.00,740.00){\line(-1,0){60.00}}
\put(40.00,740.00){\line(0,1){60.00}}
\put(40.00,800.00){\circle*{6}}
\put(100.00,800.00){\circle*{6}}
\put(100.00,740.00){\circle*{6}}
\put(40.00,740.00){\circle*{6}}
\put(30.00,820.00){\makebox(0,0)[cc]{$1$}}
\put(110.00,820.00){\makebox(0,0)[cc]{$1$}}
\put(30.00,720.00){\makebox(0,0)[cc]{$\geq 2$}}
\put(110.00,720.00){\makebox(0,0)[cc]{$\geq 2$}}
\put(70.00,770.00){\makebox(0,0)[cc]{$\tau$}}
\end{picture}
\hspace*{\fill}
\unitlength=0.4pt
\begin{picture}(80.00,80.00)(30.00,730.00)
\put(40.00,800.00){\line(1,0){60.00}}
\put(100.00,800.00){\line(0,-1){60.00}}
\put(100.00,740.00){\line(-1,0){60.00}}
\put(40.00,740.00){\line(0,1){60.00}}
\put(40.00,800.00){\circle*{6}}
\put(100.00,800.00){\circle*{6}}
\put(100.00,740.00){\circle*{6}}
\put(40.00,740.00){\circle*{6}}
\put(30.00,820.00){\makebox(0,0)[cc]{$1$}}
\put(110.00,820.00){\makebox(0,0)[cc]{$\geq 2$}}
\put(30.00,720.00){\makebox(0,0)[cc]{$\geq 2$}}
\put(110.00,720.00){\makebox(0,0)[cc]{$\geq 2$}}
\put(70.00,770.00){\makebox(0,0)[cc]{$\tau$}}
\end{picture}
\hspace*{\fill}
\unitlength=0.4pt
\begin{picture}(80.00,80.00)(30.00,730.00)
\put(40.00,800.00){\line(1,0){60.00}}
\put(100.00,800.00){\line(0,-1){60.00}}
\put(100.00,740.00){\line(-1,0){60.00}}
\put(40.00,740.00){\line(0,1){60.00}}
\put(40.00,800.00){\circle*{6}}
\put(100.00,800.00){\circle*{6}}
\put(100.00,740.00){\circle*{6}}
\put(40.00,740.00){\circle*{6}}
\put(30.00,820.00){\makebox(0,0)[cc]{$\geq 2$}}
\put(110.00,820.00){\makebox(0,0)[cc]{$\geq 2$}}
\put(30.00,720.00){\makebox(0,0)[cc]{$\geq 2$}}
\put(110.00,720.00){\makebox(0,0)[cc]{$\geq 2$}}
\put(70.00,770.00){\makebox(0,0)[cc]{$\tau$}}
\end{picture}
\hspace*{\fill}
\end{center}
For these four cases we get:
\kitem
\item
$\bigcap_{a^\ai\in\tau}\big(\span_\C K^R_{a^\ai}\big)=\tau^\bot=
\span_\C \big(\bigcap_{a^\ai\in\tau} K^R_{a^\ai}\big)$,
\item
$\bigcap_{a^\ai\in\tau}\big(\span_\C K^R_{a^\ai}\big)=
(a^3,a^4)^\bot=\tau^\bot+
\unitlength=0.3pt
\begin{picture}(80.00,60.00)(30.00,760.00)
\put(40.00,800.00){\line(1,0){60.00}}
\put(100.00,800.00){\line(0,-1){60.00}}
\put(100.00,740.00){\line(-1,0){60.00}}
\put(40.00,740.00){\line(0,1){60.00}}
\put(40.00,800.00){\circle*{6}}
\put(100.00,800.00){\circle*{6}}
\put(100.00,740.00){\circle*{6}}
\put(40.00,740.00){\circle*{6}}
\put(30.00,820.00){\makebox(0,0)[cc]{$0$}}
\put(110.00,820.00){\makebox(0,0)[cc]{$0$}}
\put(30.00,720.00){\makebox(0,0)[cc]{$1$}}
\put(110.00,720.00){\makebox(0,0)[cc]{$1$}}
\end{picture}
=\span_\C \big(\bigcap_{a^\ai\in\tau} K^R_{a^\ai}\big)$,
\item
$\bigcap_{a^\ai\in\tau}\big(\span_\C K^R_{a^\ai}\big)=
(a^4)^\bot=\tau^\bot+
\unitlength=0.3pt
\begin{picture}(80.00,90.00)(30.00,760.00)
\put(40.00,800.00){\line(1,0){60.00}}
\put(100.00,800.00){\line(0,-1){60.00}}
\put(100.00,740.00){\line(-1,0){60.00}}
\put(40.00,740.00){\line(0,1){60.00}}
\put(40.00,800.00){\circle*{6}}
\put(100.00,800.00){\circle*{6}}
\put(100.00,740.00){\circle*{6}}
\put(40.00,740.00){\circle*{6}}
\put(30.00,820.00){\makebox(0,0)[cc]{$0$}}
\put(110.00,820.00){\makebox(0,0)[cc]{$0$}}
\put(30.00,720.00){\makebox(0,0)[cc]{$1$}}
\put(110.00,720.00){\makebox(0,0)[cc]{$1$}}
\end{picture}
\,+\,
\begin{picture}(80.00,90.00)(30.00,760.00)
\put(40.00,800.00){\line(1,0){60.00}}
\put(100.00,800.00){\line(0,-1){60.00}}
\put(100.00,740.00){\line(-1,0){60.00}}
\put(40.00,740.00){\line(0,1){60.00}}
\put(40.00,800.00){\circle*{6}}
\put(100.00,800.00){\circle*{6}}
\put(100.00,740.00){\circle*{6}}
\put(40.00,740.00){\circle*{6}}
\put(30.00,820.00){\makebox(0,0)[cc]{$0$}}
\put(110.00,820.00){\makebox(0,0)[cc]{$1$}}
\put(30.00,720.00){\makebox(0,0)[cc]{$0$}}
\put(110.00,720.00){\makebox(0,0)[cc]{$1$}}
\end{picture}
=\span_\C \big(\bigcap_{a^\ai\in\tau} K^R_{a^\ai}\big)$, and
\item
$\bigcap_{a^\ai\in\tau} K^R_{a^\ai}$ contains
\unitlength=0.3pt
\begin{picture}(80.00,90.00)(30.00,760.00)
\put(40.00,800.00){\line(1,0){60.00}}
\put(100.00,800.00){\line(0,-1){60.00}}
\put(100.00,740.00){\line(-1,0){60.00}}
\put(40.00,740.00){\line(0,1){60.00}}
\put(40.00,800.00){\circle*{6}}
\put(100.00,800.00){\circle*{6}}
\put(100.00,740.00){\circle*{6}}
\put(40.00,740.00){\circle*{6}}
\put(30.00,820.00){\makebox(0,0)[cc]{$1$}}
\put(110.00,820.00){\makebox(0,0)[cc]{$1$}}
\put(30.00,720.00){\makebox(0,0)[cc]{$0$}}
\put(110.00,720.00){\makebox(0,0)[cc]{$0$}}
\end{picture}
\hspace{0.5em},\hspace{0.5em}
\begin{picture}(80.00,90.00)(30.00,760.00)
\put(40.00,800.00){\line(1,0){60.00}}
\put(100.00,800.00){\line(0,-1){60.00}}
\put(100.00,740.00){\line(-1,0){60.00}}
\put(40.00,740.00){\line(0,1){60.00}}
\put(40.00,800.00){\circle*{6}}
\put(100.00,800.00){\circle*{6}}
\put(100.00,740.00){\circle*{6}}
\put(40.00,740.00){\circle*{6}}
\put(30.00,820.00){\makebox(0,0)[cc]{$1$}}
\put(110.00,820.00){\makebox(0,0)[cc]{$0$}}
\put(30.00,720.00){\makebox(0,0)[cc]{$1$}}
\put(110.00,720.00){\makebox(0,0)[cc]{$0$}}
\end{picture}\hspace{0.5em},\hspace{0.5em}and\hspace{0.5em}
\begin{picture}(80.00,90.00)(30.00,760.00)
\put(40.00,800.00){\line(1,0){60.00}}
\put(100.00,800.00){\line(0,-1){60.00}}
\put(100.00,740.00){\line(-1,0){60.00}}
\put(40.00,740.00){\line(0,1){60.00}}
\put(40.00,800.00){\circle*{6}}
\put(100.00,800.00){\circle*{6}}
\put(100.00,740.00){\circle*{6}}
\put(40.00,740.00){\circle*{6}}
\put(30.00,820.00){\makebox(0,0)[cc]{$0$}}
\put(110.00,820.00){\makebox(0,0)[cc]{$0$}}
\put(30.00,720.00){\makebox(0,0)[cc]{$1$}}
\put(110.00,720.00){\makebox(0,0)[cc]{$1$}}
\end{picture}
\hspace{0.1em};\hspace{0.5em}hence
$\span_\C \big(\bigcap_{a^\ai\in\tau} K^R_{a^\ai}\big)=M_\C$.
\vspace{2ex}
\kenditem
So indeed $T^2_{\tau,\,\mbox{\kf loc}}(-R)=0$ for all $R$.
\par


\neu{deform-contract}
Since we have related $T^\kn_X(-R)$ to the cohomolgy groups of
$C_\kb\big(\fsigma,V^R_\kb\big)$, we are going to show the exactness of this
complex for the degrees in question. Let us begin with
a topological lemma stating the
contractibility of certain subcomplexes of polytopes.
\par

{\bf Lemma:}
{\em
Let $\kP\subseteq\R^n$ be a polytope.
For a hyperplane $H\subseteq\R^{n+1}$ and any subfan
$\CC\subseteq\big\{\tau\leq\cone(\kP)\,\big|\; \tau\subseteq H\big\}$,
we define
\[
\fcP^{H,\CC}:=\big\{\tau\leq\cone(\kP)\,\big|\; \tau\subseteq H^+\,
    \mbox{\rm  and }\, \tau\cap H\in\CC\big\}
\]
with $H^+$ denoting a closed half space corresponding to $H$.
Then, if $\kP\cap \innt(H^+)$ is non-empty, the constant cohomological
system is acyclic, i.e.\
\[
H^\kb\Big(\fcP^{H,\CC},\,\Z\Big)=0\,.
\vspace{-3ex}
\]
}
\par

{\bf Proof:}
We have to check that the corresponding polyhedral subcomplex
$\kP^{H,\CC}\subseteq\kP$ is contractible. But this is a consequence
of the following two points:
\kitem
\item[(i)]
$\kP_{H,\CC}:=\kP\cap\big[\innt( H^+) \cup |\CC|\big]\,\subseteq \kP$
is star shaped, hence contractible.
\item[(ii)]
We use the general fact that, if $Q$ is a polytope and $\widetilde{H}^+$
is a subset of the closed halfspace $H^+$ containing
$\innt(H^+)$ with $Q\not\subseteq \widetilde{H}^+$,
then $\partial Q\cap \widetilde{H}^+$ is a deformation retract of
$Q\cap \widetilde{H}^+$.
This enables us to successively get rid
of ``damaged'' $\kP$-faces contained in $\kP_{H,\CC}$. In the end we get
that $\kP^{H,\CC}$ is a deformation retract of $\kP_{H,\CC}$.
\vspace{-5ex}
\kenditem
%
%
\hspace*{\fill}$\Box$
\par


\neu{deform-geq2}
We return to the situation of \zitat{deform}{leq1} and
\zitat{deform}{loc}, i.e.\
$\kP\subseteq\R^n$ is a lattice polytope giving rise to the Gorenstein cone
$\sigma:=\cone(\kP)\subseteq N_\R=\R^{n+1}$.
\par

{\bf Proposition:}
{\em
If $R\in M$ is a degree such that $R\not\leq 1$ on $\kP$, then the
complex induced by the homological system $V^R_\kb$ is exact.
}
\par

{\bf Proof:}
The degree $R\in M$ induces a subfan
\[
\fcP^{[R\geq 1]}:=\{ \tau\leq\cone(\kP)\,|\;
\langle a^\ai,R\rangle\geq 1\;\mbox{ for every }a^\ai\in\tau\}\;\subseteq\fcP.
\]
For every $\tau\in\fcP^{[R\geq 1]}$, we write $\tilde{\tau}\leq\tau$ for the
face spanned only by those generators $a^\ai\in\tau$ satisfying 
$\langle a^\ai,R\rangle=1$. The homological system $V^R_\kb$ can more conveniently
be described  as
\[
V^R_\tau=\left\{\begin{array}{ll}
\tilde{\tau}^\bot\subseteq M_\C & \mbox{if } \tau\in \fcP^{[R\geq 1]}\\
0 & \mbox{otherwise}.
\end{array}\right.
\]
We construct a homotopy between $0$ and the identical map
$\mbox{id}:(V^R_\kb)^\ast \to (V^R_\kb)^\ast$.
Hence, denoting by $\Z\big[\fcP^{[R\geq 1]}_i\big]$ the free abelian group 
generated by the $\kn$-dimensional cones, it remains to construct a homotopy
\[
\dgARROWLENGTH=2em
\begin{diagram}
\node{\Z\big[\fcP^{[R\geq 1]}_{\kn+1}\big]}
\arrow{e,t}{\partial}
\arrow{s,l}{\mbox{\kf id}}
\node{\Z\big[\fcP^{[R\geq 1]}_{\kn}\big]}
\arrow{e,t}{\partial}
\arrow{s,l}{\mbox{\kf id}}
\arrow{sw,t}{D^\kn\hspace{-0.7em}}
\node{\Z\big[\fcP^{[R\geq 1]}_{\kn-1}\big]}
\arrow{s,l}{\mbox{\kf id}}
\arrow{sw,t}{D^{\kn-1}\hspace{-0.8em}}\\
\node{\Z\big[\fcP^{[R\geq 1]}_{\kn+1}\big]}
\arrow{e,t}{\partial}
\node{\Z\big[\fcP^{[R\geq 1]}_{\kn}\big]}
\arrow{e,t}{\partial}
\node{\Z\big[\fcP^{[R\geq 1]}_{\kn-1}\big]}
\end{diagram}
\]
such that $D^\kn(\tau)\in \Z\big[\fcP^{[R\geq 1]}_{\kn+1}\big]$ 
may be written as
$D^\kn(\tau)=\sum_v\lambda_v\phi_v$ with $\lambda_v\in\Z$ and
$\phi_v\in\fcP^{[R\geq 1]}$ such that $\tilde{\phi}\subseteq\tilde{\tau}$:\\
%
Assume that $D^{\kn-1}$ has been already constructed. If 
$\tau\in\fcP^{[R\geq 1]}$ is an $\kn$-dimensional cone, then we can
apply Lemma \zitat{deform}{contract} with $H:=[R=1]$ and 
$\CC$ being the fan consisting of $\tilde{\tau}$ and its faces.
Since $\big(\tau - D^{\kn-1}(\partial\tau)\big)\in\fcP^{H,\CC}$ and
\[
\partial\big(\tau - D^{\kn-1}(\partial\tau)\big)
\;=\;
\partial\tau - (\partial\circ D^{\kn-1})(\partial\tau)
\;=\;
(D^{\kn-2}\circ\partial)(\partial\tau)
\;= 0,
\]
there exists an element $D^\kn(\tau)\in\fcP^{H,\CC}_{\kn+1}$ such that
$\partial D^\kn(\tau)= \tau - D^{\kn-1}(\partial\tau)$.
\hfill$\Box$
\par

As a straightforward consequence of the Propositions
\zitat{deform}{leq1}, \zitat{deform}{loc}, \zitat{deform}{geq2} 
and of Theorem \zitat{van}{pyramid} 
we obtain the following
\par

{\bf Theorem:}
{\em
Assume that the two-dimensional faces of $\kP$ are either squares or 
triangles with area $1$ and $1/2$, respectively, i.e.\
$X_\kP$ is a conifold in codimension three. Then, if $R\in M$ is any
degree, we have for $\kn\leq 2$
\[
\renewcommand{\arraystretch}{1.2}
T^\kn_X(-R)=\left\{\begin{array}{@{}cl}
\kT^\kn\big(\kP\cap [R=1]\big) & \mbox{if } R\leq 1 \mbox{ on }\kP\\
0 & \mbox{otherwise}.
\end{array}\right.
\]
In particular, if $\kP$ additionally
satisfies the cleaning condition \zitat{van}{pyramid}
and contains only pyramids as three-dimensional faces, then $T^2_X=0$.
}
\par


%
%

{\small
\parbox{7cm}{
Klaus Altmann\\
Institut f\"ur Reine Mathematik\\
Humboldt-Universit\"at zu Berlin\\
Ziegelstr.~13A\\
D-10099 Berlin, Germany\\
e-mail: altmann@mathematik.hu-berlin.de}
\hfill
\parbox{6cm}{
Duco van Straten\\
Fachbereich Mathematik (17)\\
Johannes Gutenberg-Universit\"at\\
D-55099 Mainz\\
e-mail: straten@mathematik.uni-mainz.de}}


{\small
\begin{thebibliography}{BCKvS}

\bibitem[Al]{versal} Altmann, K.: The versal Deformation of an isolated toric
Gorenstein Singularity.
Invent.~math.~{\bf 128}, 443-479 (1997).

\bibitem[AH]{AH} Altmann, K., Hille, L.: 
Strong exceptional sequences provided by quivers. 
Algebras and Representation Theory {\bf 2}(1), 1-17 (1999).

\bibitem[AvS]{AvS} Altmann, K., van Straten, D.:
Quiver polytope varieties and their deformations. In preparation.

\bibitem[AS]{AQ} Altmann, K., Sletsj\o{}e, A.B.:
Andr\'{e}-Quillen cohomology of monoid algebras.
J.~of Algebra {\bf 210}, 708-718 (1998).

\bibitem[Ba]{reflexive} Batyrev, V. V.: Dual polyhedra and mirror symmetry for
Calabi-Yau hypersurfaces in toric varieties.
J.~Algebraic Geometry {\bf 3} (1994), 493-535.

\bibitem[BCKvS]{BCKvS} Batyrev, V.V.; Ciocan-Fontanine, I.; Kim, B.; 
van Straten, D.: Mirror Symmetry and Toric Degenerations of Partial 
Flag Manifolds. 
E-print math.AG/9803108. 

\bibitem[BC]{BeCh} Behnke, K., Christophersen, J.A.: Hypersurface sections and
obstructions (rational surface singularities).
Compositio Math. {\bf 77} (1991), 233-268.

\bibitem[Br]{Br} Brion, M.: The structure of the polytope algebra.
Tohoku Math.~Journal {\bf 49} (1997), 1-32.

\bibitem[Da]{Da} Danilov, V.I.: The Geometry of Toric Varieties.
Russian Math. Surveys {\bf 33}/2 (1978), 97-154.

\bibitem[GM1]{GM} Gelfand, S.I., Manin, Yu.I.: Methods of Homological 
Algebra.
Springer 1996.

\bibitem[GM2]{EMS} Gelfand, S.I., Manin, Yu.I.: Homological Algebra.
Encyclopaedia of Mathematical Sciences {\bf 38}, Algebra V; Springer 1994.

\bibitem[La]{L}  Lakshmibai, V.: Degeneration of flag varieties to 
toric varieties. 
C.R.~Acad.~Sci.~Paris {\bf 321} (1995), 1229-1234.

\bibitem[Lo]{Loday} Loday, J.-L.: Cyclic Homology.
Grundlehren der mathematischen Wissenschaften {\bf 301}, Springer-Verlag 1992.

\bibitem[St]{S} Sturmfels, B.: Gr\"obner Bases and Convex Polytopes. 
Univ.~Lect.~Notes, vol.~8, AMS, 1996. 

\end{thebibliography}
}
\end{document}